\documentclass[11pt]{amsart}
\usepackage[T1]{fontenc}
\usepackage{amssymb,mathtools,booktabs}
\usepackage[margin=1.25in,top=1in,bottom=1in]{geometry}
\usepackage{xcolor}
\usepackage{multirow}
\usepackage{tikz}
\usetikzlibrary{positioning}

\newtheorem{maintheorem}{Theorem}

\newtheorem{theorem}{Theorem}[section]
\newtheorem{proposition}[theorem]{Proposition}
\newtheorem{lemma}[theorem]{Lemma}
\newtheorem{corollary}[theorem]{Corollary}

\theoremstyle{remark}
\newtheorem{remark}[theorem]{Remark}

\numberwithin{equation}{section}

\DeclareMathOperator{\NS}{NS}

\DeclareMathOperator{\rk}{rk}

\DeclareMathOperator{\Proj}{Proj}

\DeclareMathOperator{\Trace}{Trace}

\newcommand{\ZZ}{\mathbb Z}
\newcommand{\QQ}{\mathbb Q}
\newcommand{\FF}{\mathbb F}
\newcommand{\PP}{\mathbb P}
\newcommand{\RR}{\mathbb R}

\newcommand{\et}{\mathrm{\acute et}}

\title[Pointless K3 surfaces]
{Pointless K3 surfaces over finite fields}

\author{Asher Auel}
\address{Department of Mathematics, Dartmouth College, Hanover, New Hampshire}
\email{asher.auel@dartmouth.edu}

\author{Philip Engel}
\address{Department of Mathematics, University of Illinois, Chicago}
\email{pengel@uic.edu}

\author{Kiran Kedlaya}
\address{Department of Mathematics, University of California San Diego, La Jolla, CA}
\email{kedlaya@ucsd.edu}

\author{Sair Shaikh}
\email{sairshaikh1@gmail.com}

\date{\today}


\usepackage{hyperref}
\hypersetup{
  pdftitle={Pointless K3 surfaces over finite fields},
  pdfauthor={Asher Auel and Phil Engel and Kiran Kedlaya and Sair Shaikh},
  colorlinks=true,
  linkcolor=blue,
  anchorcolor=blue,
  citecolor=blue
}

\begin{document}

\begin{abstract}
We prove that K3 surfaces with no rational points over the finite field with $q$ elements exist if and only if $q \leq 11$.
Our approach applies linear programming techniques to improve on generic sphere packing bounds for the discriminants of even rootless definite lattices.  These bounds are applied to constrain the possible zeta functions of K3 surfaces with no rational points via the Artin--Tate formula.  The remaining cases are ruled out by direct analysis of primitive embeddings of particular lattices into Niemeier lattices.
\end{abstract}

\maketitle

\section*{Introduction}

For an algebraic variety over a finite field, the number of rational points is a fundamental invariant, which has been studied for a variety of reasons. For example, algebraic curves (and recently, higher-dimensional varieties) with many rational points play an important role in the construction of error-correcting codes; for this reason among others, the study of the maximum number of points on a curve of given genus over a given finite field has developed dramatically over the past four decades (see, for example, \cite{manypoints}, \cite{faber_grantham}, \cite{castryck_voight}, \cite{kodama_top_washio}). While minimizing the number of points has fewer direct applications, it is closely related to the maximization problem, and has thus also attracted attention (see, for example, \cite{Stark}, \cite{maisner_nart}, \cite{howe_lauter_top}, \cite[Section~5]{poonen:curves}, \cite{aubry_herbaut_monaldi}, \cite{stirpe}, \cite{becker_glass}, \cite{stichtenoth}).

By contrast, while the theoretical machinery for studying the number of points on a higher-dimensional variety using the properties of its zeta function (i.e., appealing to the Weil conjectures) is well-developed, there have been comparatively fewer effective results than for curves. A natural test case is that of K3 surfaces; here we focus on the question of the existence of \emph{pointless} K3 surfaces, meaning those with no rational points over the field of definition.
One well-known example of a pointless K3 surface is the Fermat quartic surface
\[
x^4+y^4+z^4+w^4=0
\]
in $\PP^3$ over $\FF_5$.  On the other hand, the Riemann hypothesis part of the Weil conjectures implies that any K3 surface over $\FF_q$ with $q \geq 23$ has a rational point (see Proposition~\ref{prop:23}).  Here we show that this bound can be significantly improved.

\begin{maintheorem} 
\label{thm:main}
There exists a pointless K3 surface over $\FF_q$ if and only if $q \leq 11$.
\end{maintheorem}

This has an immediate consequence in higher dimensions due to the work of Bernasconi--Filipazzi \cite{BernasconiFilipazzi}; they prove a number of results on the existence of rational points on various types of threefolds over $\FF_q$ assuming that $q$ is a power of a prime $p>5$ and that $q > 19$. The former restriction comes from limitations in the minimal model program in small characteristics, but the latter is used via the hypothesis that there is no pointless K3 surface over $\FF_q$. Using Theorem~\ref{thm:main}, one can revise \cite[Proposition~3.4]{BernasconiFilipazzi} to replace the bound $q>19$ with $q > 11$; this improvement then propagates to all of the results of \cite{BernasconiFilipazzi}.

While working on this problem, we learned of the work of Brandhorst--Takada \cite{BrandhorstTakada} giving another proof of Theorem~\ref{thm:main}. We explain the structure of our proof by comparing with that of Brandhorst--Takada.
\begin{itemize}
    \item In both proofs, one starts by enumerating possible $L$-polynomials of pointless K3 surfaces. In \cite[\S 4]{BrandhorstTakada} this is implemented directly,
    whereas we use the library for this purpose described in \cite{KS} and included in SageMath.
    \item In both proofs, one continues by considering discriminants of sublattices of the Picard lattice. In \cite[\S 5]{BrandhorstTakada} this is done by analyzing a Frobenius action on the K3 lattices coming from Taelman's construction of ``Deligne modules'' for ordinary K3 surfaces \cite{Taelman}, whereas we use the Artin--Tate formula. There seems to be a close relationship between these perspectives which we have not worked out.
    \item In both proofs, one identifies that the $-1$ eigenspace of the Picard lattice for the Frobenius action is an even definite rootless lattice, and rules out many possible $L$-polynomials using discriminant bounds for such lattices. In \cite[\S 7]{BrandhorstTakada} these bounds are obtained by specializing generic sphere packing bounds, whereas we modify the proofs of these bounds to optimize them further for this particular use case; this results in an orders-of-magnitude reduction in complexity for the subsequent arguments (see Table~\ref{table:count of candidates}).
    \item In both proofs, one next rules out cases based on values whose representability by the Picard lattice is incompatible with pointlessness. In \cite[\S 8]{BrandhorstTakada}, these values are $-2,0$; we 
    treat these cases by similar methods but also use the value $2$, taking advantage of the fact that a smooth curve of genus 2 over $\FF_q$ must have a rational point if $q \geq 13$.
    \item In \cite[\S 7]{BrandhorstTakada}, the remaining cases are handled by classifying Frobenius-equivariant embeddings of the Picard lattice into the K3 lattice. We instead first rule out a few cases in which the $-1$ eigenspace of the Picard lattice is forced to admit an even integral overlattice, by introducing incompatible upper and lower bounds on the number of roots in the overlattice. We then handle the remaining cases by considering embeddings of the $-1$ eigenspace into a Niemeier lattice, but without computing any of these embeddings explicitly: we identify a linear programming constraint on such embeddings coming from the Niemeier--Venkov identity, and then construct certificates for the nonsolvability of this LP problem.
\end{itemize}

Our argument depends on several numerical computations; all necessary code and data can be found in the following GitHub repository:
\begin{center}
\texttt{https://github.com/Sair-Shaikh/pointless\_K3s}
\end{center}

\subsection*{Acknowledgments}

Thanks to Simon Brandhorst, Henry Cohn, Adam Logan, and Yuta Takada for helpful discussions, and to Brandhorst and Takada for providing an advance copy of \cite{BrandhorstTakada}.
Auel was partially supported by a McLane Family Fellowship and a Walter and Constance Burke Research Initiation Award from Dartmouth College.  Kedlaya was supported by NSF (grant DMS-2401536) and UC San Diego (Stefan E. Warschawski professorship). The authors would like to thank ICERM for hosting the special semester on \emph{Computations on K3 Surfaces and Related Varieties} in Fall 2026, where this collaboration began.

\subsection*{AI disclosure}

OpenAI GPT-6 Sol, GPT-6 Astra, and Claude Opus 5.5 were used to experiment with and debug the proof strategy. Key inputs drawn from these interactions include: some reasoning about cyclic overlattices; the use of the Niemeier--Venkov identity as the basis of an LP duality argument; generation of numerical certificates for the linear programming steps; and constructions of elliptic fibrations with specified fiber types.  During the writing process, an AI interaction made us aware of work in progress of Brandhorst--Takada, after which we reached out to them. 
All text in the final manuscript was written by the authors. 
All mathematical assertions and computational artifacts were verified by the authors in an AI-free manner.

\section{Weil bounds}

Throughout, let $q$ denote a power of a prime $p$.  For a variety $X$ defined over a finite field $\FF_q$, denote by $\overline{X}$ its base change to $\overline{\FF}_q$.  

Using the statements of the factorization of the zeta function in the Weil conjectures (proved using the Grothendieck--Lefschetz trace formula), we can extract the following general \emph{Weil bounds} on the number of points on a variety over a finite field.

\begin{theorem}[Weil bound]
\label{thm:weil_bound}
Let $X$ be a smooth, projective, geometrically irreducible variety of dimension $n$ over a finite field $\FF_q$.  Then
\[
|\#X(\FF_q) - (q^n + q^{n-1} + \dotsm + q + 1)| \leq \sum_{i=1}^{2n-1} b_i \, q^{i/2}
\]
where $b_i = \dim_{\QQ_\ell} H^i_{\et,\text{prim}}(\overline{X},\QQ_\ell)$ are the primitive $\ell$-adic Betti numbers of $X$ (for some prime $\ell$ not dividing $q$).
\end{theorem}

We recall various refinements, notably due to Serre and Stark, to the Weil bound for curves over a finite field, and their implications for curves of small genus.

Wedderburn's theorem implies that any genus zero curve, such as a smooth conic, has a rational point. We will need a slightly stronger version of this statement.
\begin{proposition} \label{lem:genus 0 curves have points}
Let $X$ be a reduced, geometrically connected (but not necessarily geometrically irreducible) curve of arithmetic genus $0$ over $\FF_q$. Then $X(\FF_q) \neq \varnothing$.
\end{proposition}
\begin{proof}
The irreducible components of $\overline{X}$ form a contractible graph of smooth rational curves.
By the Lefschetz fixed point theorem, the action of Frobenius on the topological realization of this graph has a fixed point. This corresponds to either a fixed vertex, where two rational curves meet at a unique and hence $\FF_q$-rational point; or a fixed edge, corresponding to a smooth rational curve defined over $\FF_q$ and hence having points.
\end{proof}

For genus one curves, the following is a special case of Theorem~\ref{thm:weil_bound} that predates the general formulation of the Weil conjectures.

\begin{theorem}[Schmidt--Weil]
Every smooth, projective, geometrically irreducible curve of genus one over $\FF_q$ has a rational point.    
\end{theorem}
\begin{proof}
A result of Schmidt~\cite{schmidt} can be interpreted as the statement that every smooth projective geometrically irreducible curve over a finite field has a divisor of degree 1; in genus 1, Riemann--Roch then provides a rational point. Alternatively,  Weil~\cite{weil:courbes} established Theorem~\ref{thm:weil_bound} for curves, which also implies the result.
(Note that Hasse's earlier proof~\cite{hasse:boundI,hasse:boundII} of the ``Hasse bound'' for elliptic curves  does \emph{not} imply the claim, as it requires the existence of a rational point in order to use the group law.)
\end{proof}

For higher genus curves, the Weil bound was sharpened by the following result of Serre~\cite{serre:curves_notes,serre:curves_book}.

\begin{theorem}[Serre]
\label{thm:serre}
Let $C$ be a smooth, projective, geometrically irreducible curve of genus~$g$ over $\FF_q$.  Then
\[
|\#C(\FF_q) - (q+1)| \leq g \lfloor 2 \sqrt{q} \rfloor.
\]
\end{theorem}

In genus two, we have the following further refinement.

\begin{corollary}[Serre--Stark]
\label{cor:serre-stark}
For $q$ a prime power, there exists a smooth, projective, geometrically irreducible curve $C$
of genus $2$ over $\FF_q$ with  $\#C(\FF_q) = \varnothing$ if and only if $q < 13$.  
\end{corollary}
\begin{proof}
The existence of rational points on genus two curves over $\FF_q$ for $q > 13$ follows immediately from Theorem~\ref{thm:serre}, and the case $q=13$ is due to Stark~\cite{Stark}.  Conversely, for each prime power $q < 13$, \cite[Table~4]{maisner_nart} gives an example of a pointless genus $2$ curve over~$\FF_q$.
\end{proof}

\section{Weil polynomials of K3 surfaces}
\label{sec:Weil K3}

For the rest of the paper, let $X$ be a K3 surface over $\FF_q$.
As a first step, we apply the Weil bound to obtain a nonexistence result for pointless K3 surfaces for sufficiently large $q$.
\begin{proposition}
\label{prop:23}
For any K3 surface $X$ over $\FF_q$,
\[
|\#X(\FF_q) - (q^2+1)| \leq 22q.
\]
In particular, if $q \geq 23$, then $X(\FF_q) \neq \varnothing$.
\end{proposition}
\begin{proof}
Since $b_1(X)=b_3(X)=0$ and $b_2(X)=22$, the Weil bound (Theorem~\ref{thm:weil_bound}) immediately yields the first assertion.
This in turn implies that $\#X(\FF_q) \geq q^2 - 22q + 1$, which is positive for all $q \geq 23$.
\end{proof}

To go further, we use the finer structure of the zeta function of $X$ first described by Deligne \cite{Deligne}: it has the form
\[
\frac{1}{(1-T) q^{-1} Q(qT) (1-q^2T)}
\]
for some polynomial $Q(T)$ of degree 22 over $\ZZ$ with constant term $q$ and all complex roots on the unit circle,
which we call the \emph{$L$-polynomial} of $X$.
In particular, if $\alpha_1,\dots,\alpha_{22}$ are the roots of $Q(T)$, then
\begin{equation} \label{eq:K3 point count}
\#X(\FF_q) = q^2 + 1 + q \sum_{i=1}^{22} \alpha_i;
\end{equation}
this gives back Proposition~\ref{prop:23} by applying the triangle inequality.

For $i$ a positive integer, let $\Phi_i(T)$ be the \emph{reversed} $i$-th cyclotomic polynomial, so that $\Phi_1(T) = 1-T$ rather than $T-1$. We can then factor $Q(T)$ uniquely as a product
$\prod_i \Phi_i(T)^{a_i}$
times a factor which does not vanish at any root of unity; we can describe the cyclotomic factors using the fact that the Tate conjecture is known for codimension 1 cycles on $X$ (see \cite{maulik:tate_conj,Charles,MadapusiPera,KimMadapusiPera,MPerr,IIK}) and codimension 2 cycles on $X \times X$ (see \cite{IIK}).
Namely, the group $\NS(\overline{X})$ is finite free,
the action of the Frobenius $\mathrm{Frob}_q$
on $\NS(\overline{X}) \otimes \QQ$ is semisimple, and $\det(1 - T \mathrm{Frob}_q, \NS(\overline{X})) = \prod_i \Phi_i(T)^{a_i}$.

\begin{remark}
A first consequence of the previous discussion (which does not itself depend on the Tate conjecture)
is a slight improvement on Proposition~\ref{prop:23}:
\[
|\#X(\FF_q) - (q^2+q+1)| \leq 21q.
\]
Namely, since $X$ is projective, it admits a nontrivial divisor defined over $\FF_q$, which is an eigenvector of eigenvalue $q$; hence at least one of the $\alpha_i$ is equal to $1$.
While this sharper bound still only implies the existence of rational points over $\FF_q$ for $q \geq 23$, it represents a first step in the direction of refining the Weil bound.
\end{remark}

We now give a high-level description of the proof of Theorem~\ref{thm:main}.
\begin{proof}[Proof of Theorem~\ref{thm:main}]
The ``if'' implication is fully treated in Proposition~\ref{prop:pointless k3 over F11}; we thus focus on the ``only if'' implication.
By Proposition~\ref{prop:23}, we are reduced to showing that if $q \in \{13, 16, 17, 19\}$ then 
the K3 surface $X$ over $\FF_q$ cannot be pointless.
We first use the Weil polynomial enumerator in SageMath (written by the third author; see discussion in \cite{KS}) to enumerate all choices of the $L$-polynomial $Q(T)$ consistent with the condition $\#X(\FF_q) = 0$
via \eqref{eq:K3 point count}; this list is too long to include here but can be found in our GitHub repository.
We then apply the following criteria to rule out all possible $L$-polynomials;
see Table~\ref{table:count of candidates} for a numerical summary of the progressive effects of these criteria.
\begin{itemize}
    \item Refined Artin--Tate criterion (Lemma~\ref{lem:higher Artin-Tate}).
   The polynomials meeting this criterion are separated into tables based on their subsequent disposition; see Table~\ref{table:L-polys1},
Table~\ref{table:L-polys2},
and Table~\ref{table:L-polys3}. This finishes the proof for $q=19$.

    \item Excluded norm criterion (Corollary~\ref{cor:excluded norm}).
    This excludes all polynomials in Table~\ref{table:L-polys1}
    and thereby finishes the proof for $q=16$.
    \item Cyclic overlattice criterion (Lemma~\ref{lem;cyclic overlattice constraint}). 
    This excludes all polynomials in Table~\ref{table:L-polys2} and thereby
    finishes the proof for $q=17$.
    \item Niemeier embedding criterion (Lemma~\ref{cor:Niemeier embedding constraint}). 
    This excludes all polynomials in Table~\ref{table:L-polys3} and thereby finishes the proof for $q=13$.
    \qedhere
    \end{itemize}
\end{proof}

\begin{table}[ht]

\label{table:count of candidates}
\caption{Count of numbers of possible $L$-polynomials of K3 surfaces over $\FF_q$ for $q=13, 16, 17, 19$.
In each row, the counts reflect the constraints listed in all rows up to and including the row in question.}

\begin{tabular}{c|c|c|c|c}
 & $q=13$ & $q=16$ & $q=17$ & $q=19$ \\
 \hline
Weil polynomials & 67038 & 939 & 180 & 5  \\
Refined Artin--Tate & 19 & 1 & 1 & 0 \\
Excluded norm & 12 & 0 & 1 & 0 \\
Cyclic overlattice & 8 & 0 & 0 & 0 \\
Niemeier embedding & 0 & 0 & 0 & 0
\end{tabular}

\end{table}

\section{Integral lattices}

As is typical when reasoning about K3 surfaces, most of our work consists of detailed analysis of certain integral lattices. We collect the relevant terminology here.

By an \emph{integral lattice} (or more commonly just a \emph{lattice}) we will mean a finite free abelian group $K$ equipped with a nondegenerate bilinear pairing $\langle -, - \rangle_K\colon  K\times K \to \ZZ$. The lattice $K$ is \emph{even} if $\langle v,v \rangle_K \in 2\ZZ$ for all $v \in K$ and \emph{odd} otherwise.

For $K$ even, the pairing gives rise to an injective morphism $K \to K^\vee$ of abelian groups whose cokernel $A_K$ carries an induced quadratic form $q_K\colon K \to \QQ/2\ZZ$; we call $A_K$ the \emph{discriminant group} and $(A_K, q_K)$ the \emph{discriminant form} of $K$.
The order $\#A_K$ can also be interpreted as $|\det(M)|$ where $M$ is the \emph{Gram matrix} expressing the action of the pairing on some basis of $K$.

For $n \in \ZZ$ nonzero, let $K(n)$ be the lattice with the same underlying group as $K$, but with the pairing scaled by a factor of $n$; we also use this notation for $n \in \QQ$ provided that the resulting pairing is still integral.
The isomorphism class of the real quadratic space $K \otimes \RR$ is that of $\RR^{\oplus r} \oplus \RR(-1)^{\oplus s}$ for some nonnegative integers $r,s$ with $r+s = \rk(K)$; the pair $(r,s)$ is called the \emph{signature} of $K$.
The quantity $(-1)^s \#A_K$ coincides with the \emph{discriminant} of $K$ (according to some sources but not others; we will take this as a definition).

An \emph{overlattice} of $K$ is a lattice $K'$ of the same rank equipped with an inclusion $K \hookrightarrow K'$ of quadratic spaces. When $K'/K$ is a cyclic group, we say that $K'$ is a \emph{cyclic overlattice} of $K$.
If $K$ is odd, then it is an index-2 overlattice of its maximal even sublattice.

A \emph{root} of $K$ is an element $v \in K$ with $\langle v, v \rangle_K = 2$; if no such elements exist we say that $K$ is \emph{rootless}. We define the \emph{root lattice} of $K$ to be the even sublattice $R(K)$ (possibly of smaller rank) generated by the roots of $K$; it is a direct sum of lattices corresponding to simply laced Dynkin diagrams. We correspondingly write $A_n, D_n, E_n$ to mean the root lattices defined by these root systems.

The \emph{genus} of $K$ is the (finite) set of isomorphism classes of lattices $K'$ which are everywhere locally isomorphic to $K$ (including at the infinite place, so the signatures coincide).

We introduce some lower bounds on the discriminants of even rootless lattices, which are derived using a standard technique but do not appear to be in existing literature.
We have optimized the bounds somewhat beyond what we need here in case they may be useful for other applications; see also Remark~\ref{rem:bounds from Niemeier embeddings}.
\begin{proposition} \label{prop:rootless lattice discriminant}
Let $K$ be an even rootless positive definite lattice of rank $k$ with $1 \leq k \leq 21$. Then $\#A_K \geq \delta_k$ for the value of $\delta_k$ listed in Table~\ref{table:packing}.
\end{proposition}
\begin{table}[ht]
\caption{Lower bounds for discriminants of even rootless positive definite lattices.}
\label{table:packing}
\begin{tabular}{c|cccccccccccc}
$k$ & 1 & 2 & 3 & 4 & 5 & 6 & 7 & 8 & 9 & 10 & 11 & 12 \\
\hline
$\delta_k$ & 4 & 12 & 32 & 64 & 128 & 192 & 256 & 256 & 512 & 370 & 400 & 410
\end{tabular}

\begin{tabular}{c|cccccccccccc}
$k$ & 13 & 14 & 15 & 16 & 17 & 18 & 19 & 20 & 21 & 22 & 23\\
\hline
$\delta_k$ & 406 & 370 & 318 & 256 & 194 & 134 & 86 & 50 & 26 & 12 & 4
\end{tabular}

\end{table}
\begin{proof}
For $k \leq 9$, the optimal bound for $\delta_k$ has been computed using Voronoi's algorithm for the classification of perfect lattices; see \cite{dutour} for the case $k=9$ and references therein for $k < 9$. For $10 \leq k \leq 23$,
we obtain a linear programming bound by adapting the approach of Odlyzko--Sloane \cite{OdlyzkoSloane} (for kissing numbers) and Cohn--Elkies \cite[\S 7]{CE} (for general sphere packings), following the notation of the latter.
Let $f \colon \RR^k \to \RR$ be a nonzero Schwartz function such that
\begin{itemize}
\item $\hat{f}(t) \geq 0$ for all $t \in \RR^k$ (and hence $f(0) > 0$);
\item $f(x) \leq 0$ whenever $|x|^2 \in \{4,6,8,\dots\}$.
\end{itemize}
Then for any isometric embedding $K \hookrightarrow \RR^k$,
by Poisson summation we have
\[
f(0) \geq \sum_{x \in K} f(x) = \frac{1}{\sqrt{\#A_K}} \sum_{y \in K^\vee} \hat{f}(y) \geq \frac{1}{\sqrt{\#A_K}}  \hat{f}(0)
\]
and hence the lower bound
\[
\#A_K \geq \left( \frac{\hat{f}(0)}{f(0)} \right)^2.
\]
In addition to rounding up to the nearest integer, we may also use the fact that $\#A_K$ is even when $k$ is odd.

To construct test functions, we introduce the generalized Laguerre polynomials
\[
L_n^{(\alpha)}(u) = \sum_{j=0}^n \binom{n+\alpha}{n-j} \frac{(-x)^j}{j!}\qquad (n =0,1,\dots);
\]
for $\alpha := \tfrac{k}{2} - 1$, $L_n^{(\alpha)}(2\pi |x|^2) e^{-\pi|x|^2}$ is an eigenfunction for the Fourier transform with eigenvalue $(-1)^n$ .
We may thus take test functions $f$ of the form
\[
f(x) = P(2 \pi |x|^2) e^{-\pi|x|^2},
\qquad
\hat{f}(y) = Q(2 \pi |y|^2) e^{-\pi|y|^2}
\]
where $P$ and $Q$ are polynomials of the form
\[
Q(u) = \sum_{n=0}^\infty q_n L_j^{(\alpha)}(u), \qquad
P(t) = \sum_{n=0}^\infty (-1)^n q_n L_j^{(\alpha)}(t)
\]
for some $q_n \in \QQ$ (not all zero) for which
\begin{itemize}
\item $Q(u) \geq 0$ for all $u \geq 0$;
\item $P(t) \leq 0$ whenever $t = 4\pi m$ for some integer $m \geq 2$.
\end{itemize}
These conditions can be checked by verifying that
\begin{itemize}
     \item $Q(0) > 0$ and every real root of $Q$ is negative;
     \item $P$ has negative leading coefficient and $P(4 \pi m) < 0$ for all integers $m > 2$
     such that $4 \pi m$ is less than the largest real root of $Q$.
\end{itemize}
We provide suitable numerical values, and a SageMath script to verify their accuracy, in our GitHub repository. 
\end{proof}

\begin{remark}
It is evident from our proof of Proposition~\ref{prop:rootless lattice discriminant} that the bound is optimal for $k \leq 9$. It is also optimal for $k=16$, as witnessed by the Barnes--Wall lattice; it may be worth investigating whether the test function in this case is approximating a ``magic function'' with special properties, as in the work of Viazovska \cite{Viazovska}.

A related point is that when $K$ embeds into the Leech lattice, its orthogonal complement $K'$ is also even and rootless and $A_K \cong A_{K'(-1)}$. 
For $k \leq 9$, this construction applies to the optimal lattice and gives the best known example in dimension $24-k$; this shows that the bounds for $k = 22, 23$ are optimal.
See also Remark~\ref{rem:bounds from Niemeier embeddings}.
\end{remark}

The same construction can be used to count roots based on the discriminant.
\begin{lemma} \label{lem:disciminant root count}
Let $K$ be an even positive definite lattice of rank $16$  (resp. $17$) 
with $\#A_K \leq 96$ (resp. $\#A_K \leq 77$).
Then $K$ has at least $22$ roots.
\end{lemma}
\begin{proof}
Since roots come in inverse pairs, it is equivalent to show that if $K$ has at most 20 roots,
then $\#A_K \geq 97$ (resp. $\#A_K \geq 78$).
This is obtained by modifying the proof of Proposition~\ref{prop:rootless lattice discriminant}: with notation as therein, provided that $f(2) > 0$ we have
\[
f(0) + 20 f(2) \geq \sum_{x \in K} f(x) = \frac{1}{\sqrt{\#A_K}}
\sum_{y \in K^\vee} \hat{f}(y) \geq \frac{1}{\sqrt{\#A_K}} \hat{f}(0)
\]
and hence
\[
\#A_K \geq \left( \frac{\hat{f}(0)/f(0)}{1 + 20 f(2)/f(0)} \right)^2.
\]
We may thus read the claimed bounds bound off from the certificates provided for $k=16, 17$.
\end{proof}

\section{Application of the Artin-Tate formula}

It was observed by Elsenhans--Jahnel \cite{EJ} that the Artin--Tate formula for a K3 surface imposes nontrivial constraints on the Weil polynomial of a K3 surface over $\FF_q$ when the Picard rank does not jump from $\FF_q$ to $\FF_{q^2}$. This constraint is useless in the present context because pointless K3 surfaces necessarily have much higher Picard rank over $\FF_{q^2}$ than over $\FF_q$; however, the first and fourth authors \cite{AuelShaikh} have introduced an adaptation of the Elsenhans--Jahnel construction that gives a meaningful constraint when the Picard rank jumps from $\FF_q$ to $\FF_{q^2}$.

We first fix some additional running notation for the rest of the paper,
augmenting that introduced in \S\ref{sec:Weil K3}.
Define the \emph{Picard lattice} $P$ of the K3 surface $X$ over $\FF_q$ to be the group $\NS(\overline{X})$ equipped 
with the intersection pairing. 
By the Riemann--Roch and Hodge index theorems, $P$ is an even lattice of signature $(1, \rk(P)-1)$. 
Let $P_+, P_-$ be the eigenspaces of $P$ of eigenvalue $+1, -1$ respectively, and let  $P_{\pm}$ be the sublattice of $P$ generated by $P_+$ and $P_-$.
The following result goes back to Nikulin~\cite{Nikulin:finite_groups_short,Nikulin:finite_groups_long} and plays a pivotal role here.

\begin{lemma} \label{lemma:negative eigenspace rootless}
The lattice $P_+$ is even of signature $(1, a_1-1)$.
The lattice $P_-$ is even of signature $(0, -a_2)$ and $P_-(-1)$ is rootless.
\end{lemma}
\begin{proof}
The lattice $P_-$ is even because $P$ is.
By identifying $P_+$ with the Picard lattice of $X$, we obtain an ample class $v \in P_+$. Since $P_-$ is orthogonal to $P_+$, the signature of $P$ now forces $P_-$ to be negative definite.

Suppose by way of contradiction that there exists $w \in P_-$ with norm $-2$. By Riemann-Roch either $w$ or $-w$ is effective; without loss of generality suppose that $w$ is effective. Then $v \cdot w$ is both positive (because $v$ is ample) and zero (because $P_-$ is orthogonal to $P_+$), a contradiction.
\end{proof}

The discriminants of the lattices $P_+, P_-$ are controlled by the $L$-polynomial via the following result.

\begin{proposition}[Refined Artin--Tate formula] \label{prop:refined Artin-Tate}
Factor the $L$-polynomial $Q(T)$ of $X$ as $(1-T)^{a_1} (1+T)^{a_2} R(T)$ with $R(1) R(-1) \neq 0$
(but $R$ may have zeroes at other roots of unity).
Then for some integer $d \in \{0, \dots, \min\{a_1,a_2\}\}$, 
$\#A_{P_+}$ divides $2^d R(1)$ and
$\#A_{P_-}$ divides $2^d R(-1)$.
In addition, the integers
\[
\frac{2^{a_2} R(1)}{ \#A_{P^+}}, \frac{2^{a_2} R(-1) }{\#A_{P_-}}
\]
are perfect squares.
\end{proposition}
\begin{proof}
See \cite{AuelShaikh}.
\end{proof}

Combining with our discriminant bounds for even rootless lattices yields the following, which is a generalization of a result in \cite{AuelShaikh}.

\begin{corollary} \label{cor:refined artin-tate bound}
For any K3 surface $X$ over $\FF_q$,
$2^{\min\{a_1,a_2\}} R(-1) \geq \delta_{a_2}$
for $\delta_k$ as in Proposition~\ref{prop:rootless lattice discriminant}.
\end{corollary}
\begin{proof}
    Let $K$ be the negation of $P_-$, which by Lemma~\ref{lemma:negative eigenspace rootless} is even, positive definite, and rootless of rank $k = a_2$.
    We may therefore apply Proposition~\ref{prop:refined Artin-Tate} and Proposition~\ref{prop:rootless lattice discriminant} to obtain
    \[
    2^{\min\{a_1, a_2\}} R(-1) \geq \#A_{P_-} = \#A_K \geq \delta_{a_2}.
    \qedhere
    \]
\end{proof}

\begin{remark}
As noted earlier, the bound of Proposition~\ref{prop:rootless lattice discriminant} is optimal for $k=16$, achieved by the Barnes--Wall lattice. While that candidate for $P_-$ will be killed off relatively early in our narrative, its ghost will haunt us thereafter: the existence of such a dense lattice obstructs various strategies for ruling out even rootless lattices of slightly larger covolume.
\end{remark}

\begin{lemma} \label{lem:higher Artin-Tate}
If $q \geq 13$ and $X(\FF_q) = \varnothing$, then the $L$-polynomial of $X$ appears in one of Table~\ref{table:L-polys1}, 
Table~\ref{table:L-polys2}, 
or Table~\ref{table:L-polys3}. 
\end{lemma}
\begin{proof}
Again by Proposition~\ref{prop:23}, we need only consider $q \in \{13, 16, 17, 19\}$.
We filter the candidate $L$-polynomials found by SageMath 
using Corollary~\ref{cor:refined artin-tate bound}. In each case, $a_1 \leq a_2$ and 
$2^{a_1-1} R(-1) < \delta_{a_2}$, so 
\begin{equation} \label{eq:discriminant constraints}
\#A_{P_+} | 2^{a_1} R(1), \qquad
\sqrt{\frac{2^{a_1} R(1)}{\#A_{P_+}}} \in \ZZ, \qquad \#A_{P_-} = 2^{a_1} R(-1).
\end{equation}
This gives the list in question.
\end{proof}

\begin{table}[ht]
\caption{Candidate $L$-polynomials of pointless K3 surfaces that meet the refined Artin--Tate criterion but are ruled out because $P_+$ represents one of $-2,0,2$.
As in Proposition~\ref{prop:refined Artin-Tate}, define $R(T)$ by writing the $L$-polynomial as $(1-T)^{a_1} (1+T)^{a_2} R(T)$ with $R(1), R(-1) \neq 0$.}
\label{table:L-polys1}
\begin{tabular}{cccc}
$L$-polynomial & $R(1)$ & $\#A_{P_+}$ &  $ R(-1)$  \\
\hline
$\Phi_1 \Phi_2^{19} (16T^2 - 31T + 16)$ & 1 & 2 & 63 \\
$\Phi_1 \Phi_2^{17} \Phi_6 (13T^2 - 25T + 13)$ & 1 & 2 & 153 \\
$\Phi_1^3 \Phi_2^{17} (13T^2 - 12T + 13)$ & 14 & 28, 112 & 38  \\
$\Phi_1^2 \Phi_2^{16} (13T^4 - 12T^3 + 23T^2 - 12T + 13)$ & 25 & 1, 4, 25, 100 & 73  \\
$\Phi_1^2 \Phi_2^{16} (13T^4 - 12T^3 + 21T^2 - 12T + 13)$ & 23 & 92 & 71  \\
$\Phi_1^2 \Phi_2^{16} (13T^4 - 12T^3 + 20T^2 - 12T + 13)$ & 22 & 88 & 70 \\
$\Phi_1^2 \Phi_2^{16} (13T^4 - 12T^3 + 17T^2 - 12T + 13)$ & 19 & 76 & 67  \\
$\Phi_1^2 \Phi_2^{16} (13T^4 - 12T^3 + 14T^2 - 12T + 13)$ & 16 & 1, 4, 16, 64 & 64  \\
\end{tabular}
\end{table}

\begin{table}[ht]
\caption{Candidate $L$-polynomials of pointless K3 surfaces that meet the refined Artin--Tate and excluded norms criteria but are ruled out using cyclic overlattices.
The lattice $P_+$ is specified by its upper triangular Gram matrix: $(a,b,c)$ corresponds to the matrix $\scriptsize \begin{pmatrix} a & b \\ b & c \end{pmatrix}$. Other notation is as in Table~\ref{table:L-polys1}.}
\label{table:L-polys2}
\begin{tabular}{ccccc}
$L$-polynomial & $R(1)$ & $\#A_{P_+}$ &  $ R(-1)$ & $P_+$ \\
\hline
\multirow{2}{*}{$\Phi_1 \Phi_2^{19} (17T^2 - 16T + 17)$} & \multirow{2}{*}{18} & 4 & \multirow{2}{*}{50} & $(4)$ \\
& & 36 & & $(36)$ \\
\hline
$\Phi_1^2 \Phi_2^{16} \Phi_4 (13T^2 - 12T + 13)$ & 28 & 112 & 76 & $(4,12; 8)$ \\
\hline
\multirow{2}{*}{$\Phi_1^2 \Phi_2^{16} (13T^4 - 12T^3 + 22T^2 - 12T + 13)$} & \multirow{2}{*}{24} & \multirow{2}{*}{96} & \multirow{2}{*}{72} & $(4,8;-8)$ \\
& & & & $(8,4;-10)$ \\
\hline
\multirow{3}{*}{$\Phi_1^2 \Phi_2^{16} (13T^4 - 12T^3 + 18T^2 - 12T + 13)$} & \multirow{3}{*}{20} & 20 & \multirow{3}{*}{68} & $(4,2;-4)$ \\
& & 80 & & $(4, 8; -4)$ \\
& & 80 & & $(8,4;-8)$ \\
\hline
$\Phi_1^2 \Phi_2^{16} (13T^4 - 12T^3 + 15T^2 - 12T + 13)$ & 17 & 68 & 65 & $(8, 6; -4)$ \\
\end{tabular}
\end{table}

\begin{table}[ht]
\caption{Candidate $L$-polynomials of pointless K3 surfaces that meet the refined Artin--Tate, excluded norms, and cyclic overlattice criteria but are ruled out using embedding into Niemeier lattices.
``Genera'' counts genera of $P_-$ according to Lemma~\ref{L:genus pminus}.
Other notation is as in  Table~\ref{table:L-polys2}.}
\label{table:L-polys3}
\begin{tabular}{ccccccc}
$L$-polynomial & $R(1)$ & $\#A_{P_+}$ &  $ R(-1)$ & $P_+$ & Genera\\
\hline
$\Phi_1 \Phi_2^{17} (13T^4 - 38T^3 + 53T^2 - 38T + 13)$ & 3 & 6 & 155 & $(6)$ & 2 \\
$\Phi_1 \Phi_2^{17} (13T^4 - 38T^3 + 52T^2 - 38T + 13)$ & 2 & 4 & 154 & $(4)$ & 4 \\
$\Phi_1^2 \Phi_2^{16} (13T^4 - 12T^3 + 28T^2 - 12T + 13)$ & 30 & 120 & 78 & $(4, 0; -30)$ & 4 \\
$\Phi_1^2 \Phi_2^{16} (13T^4 - 12T^3 + 27T^2 - 12T + 13)$ & 29 & 116 & 77 & $(4, 2; -28)$ & 2 \\
$\Phi_1^2 \Phi_2^{16} \Phi_6 (13T^2 + T + 13)$ & 27 & 108 & 75 & $(6, 6; -12)$ & 4 \\
$\Phi_1^2 \Phi_2^{16} (13T^4 - 12T^3 + 24T^2 - 12T + 13)$ & 26 & 104 & 74 & $(4, 0; -26)$ & 2 \\
$\Phi_1^2 \Phi_2^{16} (13T^4 - 12T^3 + 19T^2 - 12T + 13)$ & 21 & 84 & 69 & $(4, 2; -20)$ & 2 \\
$\Phi_1^2 \Phi_2^{16} (13T^4 - 12T^3 + 16T^2 - 12T + 13)$ & 18 & 72 & 66 & $(6, 0; -12)$ & 4 \\
\end{tabular}
\end{table}

\section{The excluded norm criterion}
\label{sec:excluded norms}

We next show that pointlessness imposes some strong restrictions on the indefinite lattice $P_+ \cong \NS(X)$. 
We use the following statement that adapts a standard fact about K3 surfaces over algebraically closed fields \cite[Corollary~8.2.11]{Huybrechts}.
\begin{lemma} \label{lem:nef cone}
Let $R_X$ be the invariant part of the Weyl group of $\overline{X}$ for the action of the absolute Galois group $G_{\FF_q}$ by conjugation. 
\begin{enumerate}
    \item [(a)]
The nef positive cone of $X$ is a fundamental domain for the action of $R_X$ on the positive cone, and this action is faithful.
\item[(b)]
Every wall of the nef positive cone is the orthogonal complement to the sum over a Galois orbit $I$ of $(-2)$-curves on $\overline{X}$ of one of the following forms. (Note that every $(-2)$-curve on $\overline{X}$ is smooth of genus $0$.)
\begin{enumerate}
    \item[(i)] Any two distinct curves in $I$ are pairwise disjoint.
    \item[(ii)] There is a partition of $I$ into pairs such that the two members of each pair have intersection number $-1$ and the unions of any two distinct pairs are disjoint.
\end{enumerate}
\end{enumerate}
\end{lemma}
\begin{proof}
    See \cite[Proposition~3.6, Proposition~3.7, Remark 3.9]{BLvL}. (Note that although some results of \cite{BLvL} require characteristic not equal to 2, this does not.)
\end{proof}

\begin{proposition} \label{prop:square 0 or 2}
Suppose that $\#X(\FF_q) = \varnothing$ and $q \geq 13$. Then $\NS(X)$ admits no classes of square $-2$, $0$, or $2$.
\end{proposition}

\begin{proof}
Assume first that $\NS(X)$ admits a square $0$ class.  Then $X$ admits an elliptic fibration $X \to \PP^1$. (Here we are using that if $p \leq 3$ we cannot have a quasi-elliptic fibration, as then each fiber has a rational point by Proposition~\ref{lem:genus 0 curves have points}.)
If any of the fibers above the rational points $\PP^1(\FF_q)$ are geometrically irreducible (i.e., smooth or of type $I_1$ or $II$), then they have rational points.  Indeed, for smooth fibers this follows from Proposition~\ref{lem:genus 0 curves have points}, while for a fiber of type $I_1$ or $II$, the unique singular point of the fiber will be a rational point.  Thus we are reduced to the case when all the fibers above $\PP^1(\FF_q)$ are singular of type other than $I_1$ or $II$.  In this case, each singular fiber has multiplicity at least 2, so the maximum number of singular fibers is 12.
Since $q \geq 13$, $\#\PP^1(\FF_q) > 12$ and hence there will be a smooth fiber over a rational point, giving a contradiction. 

Next, assume that $\NS(X)$ admits a square $-2$ class.
By Riemann--Roch, either this class or its negation is effective, and thus is a positive root of $\NS(\overline{X})$.
The corresponding reflection $r$ is invariant under the Galois action and so belongs to $R_X$.
The action of $r$ on $\NS(\overline{X}) \otimes \RR$ is reflection through the hyperplane $\overline{H}$ normal to $r$.
By Lemma~\ref{lem:nef cone}(a), the action of $r$ on $\NS(X) \otimes \RR$
is also faithful, so the fixed subspace $H := \overline{H} \cap (\NS(X) \otimes \RR)$ must again have codimension 1. By Lemma~\ref{lem:nef cone}(a) again,
the nef positive cone of $X$ is a fundamental domain for the action of $R_X$,
so $H$ must abut at least one member of the $R_X$-orbit of the nef positive cone. By Lemma~\ref{lem:nef cone}(b), this means that our original positive root must be $R_X$-conjugate to the sum over a Galois orbit $I$ of $(-2)$-curves of type (i) or (ii).
For that sum to be a $(-2)$-class, $I$ must consist of a single $(-2)$-curve in case (i) or a single pair in case (ii); in both cases, by Proposition~\ref{lem:genus 0 curves have points} the resulting configuration contains an $\FF_q$-rational point, a contradiction.

Finally assume that $\NS(X)$ admits a square $2$ class. By Lemma~\ref{lem:nef cone}, the $R_X$-orbit of the chosen class contains a square $2$ class which is also nef (and obviously big).
This class defines a possibly singular model of $X$ as a double cover $\pi \colon X' \to \PP^2$; in characteristic 2, $\pi$ cannot be purely inseparable as this would imply that $\#X(\FF_q) = \#\PP^2(\FF_q) > 0$. While the map $X \to X'$ need not be an isomorphism, none of the singular geometric points of $X'$ can be $\FF_q$-rational, as otherwise Proposition~\ref{lem:genus 0 curves have points} would imply that the corresponding exceptional divisor contains an $\FF_q$-rational point. In particular, $X'(\FF_q)$ is also empty.

Let $\pi': X'' \to \PP^2$ be the quadratic twist of $\pi$.
Write $X''$ as the hypersurface
\[
w^2 + A_3(x,y,z) w = B_6(x,y,z)
\]
in the weighted projective space $\PP(3,1,1,1) = \Proj \FF_q[w,x,y,z]$ where $A_3$ and $B_6$ are homogeneous of degrees 3 and 6 respectively. 
The restriction of $\pi'$ to each $\FF_q$-rational line $\ell$ has $2q+2$ points, which by Corollary~\ref{cor:serre-stark} is too many for the cover to be geometrically irreducible.
Hence the cover splits; that is, there exists a cubic form $C_\ell(x,y,z)$ (depending on $\ell$) such that $C_\ell^2 + A_3 C_\ell - B_6$ vanishes identically on $\ell$ (not just at $\FF_q$-rational points).

We apply this condition to the line through the fixed point $[1:0:0]$ and the variable point $[0:u:v]$
where $u,v \in \FF_q$ are not both zero. We may normalize the resulting cubic form $C_{u,v}$ to have the form
\[
C_{u,v} = cx^3 + x^2 D_{1,u,v}(y,z) + x D_{2,u,v}(y,z) + D_{3,u,v}(y,z)
\]
for some fixed nonzero $c \in \FF_q$ and some forms $D_{i,u,v}(y,z)$ of degree $i$. 
By applying Hensel's lemma over the Laurent series ring $\FF_q[y,z]((x^{-1}))$,
we obtain a \emph{unique} choice of forms $D_i(y,z)$ of degree $i$ for $i=1,2,3$ for which the quantity
\[
C = cx^3 + x^2 D_1(y,z) + x D_2(y,z) + D_3(y,z)
\]
satisfies
\[
\deg_x (C^2 + A_3 C - B_6) < 3.
\]
Since this application of Hensel's lemma commutes with specialization to $u,v \in \FF_q$,
we deduce that $D_{i,u,v} = D_i$ for all $i,u,v$.
This means that the coefficients of $x^2, x^1, x^0$ in $C^2 + A_3 C - B_6$ vanish at all $y,z \in \FF_q$; since each of these is a form in $y,z$ of degree at most $6 < q+1$, these forms must vanish identically.
We deduce that $X''$ is disconnected, a contradiction.
\end{proof}

\begin{remark}
    The proof of Proposition~\ref{prop:square 0 or 2} suggests that one can construct a pointless K3 surface over $\FF_{11}$ by exhibiting an elliptic fibration
    with 12 singular fibers of type $I_2$ (and no others), each lying over a point of $\mathbb{P}^1(\FF_{11})$ and having no rational points.
    This is in fact possible, and a similar construction yields pointless K3 surfaces 
    over $\FF_q$ for all $q \leq 11$;
    see Proposition~\ref{prop:pointless k3 over F11}.
\end{remark}

Proposition~\ref{prop:square 0 or 2} has the following consequences for our candidate $L$-polynomials.
\begin{corollary} \label{cor:excluded norm}
There is no K3 surface over a finite field with $L$-polynomial appearing in Table~\ref{table:L-polys1}.
Moreover, for a K3 surface $X$ with $L$-polynomial appearing in Table~\ref{table:L-polys2} or Table~\ref{table:L-polys3}, $P_+$ must have one of the forms listed therein.
\end{corollary}
\begin{proof}
We break into cases based on the value of $a_1$.
\begin{itemize}
    \item When $a_1 = 1$, we have $\#A_{P_+} \in \{2,4,6,36\}$, but we can now argue that $\#A_{P_+} \neq 2$.
    \item 
    When $a_1 = 2$, we can reason about $P_+$ using binary quadratic forms: to the Gram matrix $\scriptsize \begin{pmatrix} 2a & b \\ b & 2c \end{pmatrix}$ we associate the norm form $2ax^2 + 2bxy + 2cy^2$ of discriminant $4ac-b^2 = -\# A_{P_+}$.
    In this notation, isomorphisms of lattices correspond to (improper) equivalences of indefinite binary quadratic forms.
    We may now extract information\footnote{We note a typo in the row labeled $-88$ in \cite[Table~15.2b]{CS}: the last entry $\{4^8 6\}$ should not be there.} from \cite[\S  15.3, Theorem 1]{CS}:
    $P_+$ must correspond to an entry of \cite[Tables~15.2a, 15.2b]{CS} for the row labeled by $-\#A_{P_+}$
    in which each subscripted entry is an even integer greater than 2.
    In addition, when $\#A_{P_+}$ is divisible by 4, the integral lattice $P_+(\frac{1}{2})$ corresponds to an entry in the row labeled by $-\frac{1}{4} \#A_{P_+}$.
    This rules out $\#A_{P_+} \in \{1, 4, 5, 8, 12, 16, 17, 24, 28, 64, 76, 88, 92, 100 \}$.

    \item There is a unique case with $a_1 = 3$, in which case $\#A_{P_+} \in \{28, 112\}$. 
    We query LMFDB for all lattices $L$ of signature $(2,1)$ 
    with $\#A_L \in \{28, 112\}$ (there are 2 of the former and 5 of the latter, all of class number 1)
    and observe that each such lattice represents one of $-2, 0, 2$.
    \qedhere
\end{itemize}

\end{proof}

In the remaining cases, 
the identification of $P_+$ will provide a constraint on $P_-$ via the following observation.

\begin{lemma} \label{lem:isom quadratic spaces}
Suppose that the $L$-polynomial of $X$ appears in 
Table~\ref{table:L-polys2} or 
Table~\ref{table:L-polys3}.
Then there exists an inclusion $A_{P_+}[2] \to A_{P_-(-1)}[2]$ of quadratic spaces.
\end{lemma}
\begin{proof}
The cokernel of the map $P_+ \oplus P_- \to P_{\pm}$ is an elementary 2-group which injects into both $A_{P_+}$ and $A_{P_-}$. In all indicated cases this group has order $2^{a_1}$, so the image in $A_{P_+}$ must coincide with the Sylow 2-subgroup of $A_{P_+}$.
\end{proof}

\section{Cyclic overlattices}
\label{sec:cyclic}

Although it is not strictly necessary to do so (as all remaining cases can be handled in the framework of embeddings into Niemeier lattices described in the next section),
we find it enlightening to eliminate a few more cases in which $P_-$ is forced to admit an even overlattice. While the latter is not guaranteed to be rootless, we can still exploit tension between lower and bounds on the number of roots. We first illustrate with a single example.

\begin{lemma}[Elkies] 
\label{lem:root count elkies}
Let $K$ be a unimodular positive definite lattice of rank $k$ which does not represent $1$.
Then $K$ has at least $2k(23-k)$ roots.
\end{lemma}
\begin{proof}
See \cite[Theorem 1]{Elkies}.
\end{proof}

\begin{remark} \label{remark:elkies bound}
The proof of \cite[Theorem~1]{Elkies} does not follow the approach of Lemma~\ref{lem:disciminant root count}; rather, it uses theta functions and modular forms. It may be feasible to use a similar approach in other cases to obtain some improvements on Lemma~\ref{lem:disciminant root count} for specific values of the discriminant.
\end{remark}

\begin{lemma} \label{lem:cyclic root count}
Let $K$ be a positive definite lattice of rank $k$ admitting a rootless sublattice, the quotient by which is cyclic of order $m$. Then every irreducible component of the root system of $K$ has Coxeter number at most $m$.
In particular, $K$ has at most $km$ roots.
\end{lemma}
\begin{proof}
Let $R_1, R_2, \dots$ be the components of the root lattice $R(K)$,
and let $h(R_i)$ be the Coxeter number of $R$.
The map $R_i \to K \to \ZZ/m\ZZ$ has nonzero image on all of the roots of $K$;
by a theorem of Kamiya--Takemura--Terao \cite[Corollary~3.4]{KTT},
this is only possible if $h(R_i) \leq m$.
(For $m=2$ one can also prove this more directly.)
The total number of roots is then
\[
\sum_i \rk(R_i) h(R_i) \leq m \sum_i \rk(R_i) = m k
\]
as desired.
\end{proof}

\begin{lemma} \label{lem;cyclic overlattice constraint1}
There is no K3 surface over a finite field with $L$-polynomial $\Phi_1 \Phi_2^{19} (17T^2 - 16T + 17)$.
\end{lemma}
\begin{proof}
Suppose to the contrary that $X$ is such a K3 surface.
The isomorphisms $A_{P_-} \cong A_{P_-}[2] \oplus A_{P_-}[5]$ and $A_{P_-}[2] \cong A_{P_+}[2]$
of quadratic spaces (the latter from Lemma~\ref{lem:isom quadratic spaces})
allow us to use a global constraint (the \emph{oddity condition} in the terminology of \cite[Chapter~15, \S 7.7]{CS}) to deduce that $A_{P_-}[5]$ must contain an isotropic vector.
We thus have inclusions $P_-(-1) \subset L \subset U$ of lattices with $[L:P_-(-1)] = 5, [U:L] = 2$;
in particular, $U$ is unimodular and $L$ is its even sublattice.
 Write $U \cong \ZZ^a \oplus U_0$ where $U_0$ does not represent $1$; then Lemma~\ref{lem:root count elkies} implies that $U_0$ has at least $2(19-a)(4+a)$ roots. Combined with the $2a(a-1)$ roots of $\ZZ^a$, this means that the number of roots of $U$ (and hence of $L$, since any root is even) is at least 
\[
2a(a-1) + 2(19-a)(4+a) = 152+28a \geq 152,
\]
whereas Lemma~\ref{lem:cyclic root count} gives a contradictory upper bound of $5 \times 19 = 95$.
\end{proof}

We next adapt this paradigm to handle a few cyclic overlattices of index $2$.

\begin{lemma} \label{lem:roots upper bound}
Let $K$ be an even positive definite lattice with $\rk(K) + \ell_2(A_K) \leq 19$ containing a rootless sublattice of index $2$. Then $K$ has at most $20$ roots.
(In particular this applies when $\rk(K) + v_2(\#A_K) \leq 19$.)
\end{lemma}
\begin{proof}
By Lemma~\ref{lem:cyclic root count}, the root lattice $R(K)$ 
has the form $A_1^a$ for some nonnegative integer $a$,
and $K$ therefore has $2a$ roots. It thus suffices to check that $a \leq 10$.

Let $S$ be the saturation of $R$ in $K$; we have a natural inclusion
\[
S \hookrightarrow R^\vee \cong \tfrac{1}{2} R.
\]
Via the induced inclusion
\[
S/R \hookrightarrow R^\vee/R \cong \#A_R \cong (\ZZ/2\ZZ)^a,
\]
we may view $S/R$ as a binary code $C$ of length $a$; for $r$ the dimension of $C$,
\[
A_S \cong (\ZZ/2\ZZ)^{a-2r}.
\]
By hypothesis, we have
\[
a-2r = \ell_2(A_S) \leq \ell_2(A_K) + \rk(S^\perp) \leq \ell_2(A_K) + \rk(K) - \rk(S) \leq 19 - a,
\]
or in other words $r \geq a - 9$.

In the other direction, by again viewing $S$ as a subgroup of $\tfrac{1}{2} R$, we see that every codeword of $C$ of weight $w$ lifts to an element of $\tfrac{1}{2} R$ of square $w/2$. This implies firstly that the code $C$ is \emph{doubly even}, meaning that each codeword has Hamming weight divisible by 4; and secondly that its minimum weight is at least 8, as a codeword of weight 4 would lift to a root in $S$
whereas every any root projects to zero in $S/R$ by definition. We can now apply the Griesmer bound \cite{Griesmer} (see also \cite[Chapter~17, Theorem~24]{MacWilliamsSloane}): a binary code of dimension $r$ and minimum weight $d$ has length at least $\sum_{i=0}^{r-1} \left\lceil \frac{d}{2^i} \right\rceil$.
For the code $C$, this yields
\[
a \geq \sum_{i=0}^{r-1} \left\lceil \frac{8}{2^i} \right\rceil = \begin{cases} 12 & r=2 \\ 14 & r = 3 \\ r+11 & r \geq 4. \end{cases}
\]
Since the right hand side is at least $r+10 > a$ for $r \geq 2$, we deduce that $r \leq 1$ and so $a \leq 10$ as desired.
\end{proof}

\begin{lemma} \label{lem;cyclic overlattice constraint}
There is no K3 surface over a finite field with $L$-polynomial appearing in Table~\ref{table:L-polys2}.
Moreover, if $X$ has $L$-polynomial appearing in Table~\ref{table:L-polys3}, then the quadratic form on $A_{P_-}[2]$ is anisotropic.
\end{lemma}
\begin{proof}
As the $L$-polynomial $\Phi_1 \Phi_2^{19} (17T^2 - 16T + 17)$ is excluded by Lemma~\ref{lem;cyclic overlattice constraint1}, we may focus on the remaining cases in Table~\ref{table:L-polys2} and Table~\ref{table:L-polys3}, for which $a_2 \in \{16, 17\}$.
Recall that by Corollary~\ref{cor:excluded norm}, $P_+$ has the form listed in the tables.

If $A_{P_-}[2]$ is isotropic, then $P_-$ admits
an even index-$2$ overlattice; we may then apply Lemma~\ref{lem:disciminant root count} and 
Lemma~\ref{lem:roots upper bound} (after noting that each listed case satisfies the hypotheses of both lemmas) to obtain a contradiction. This immediately implies the assertion about Table~\ref{table:L-polys3}.
As for Table~\ref{table:L-polys2}, in each case either $P_+$ primitively represents a multiple of 8 or $\#A_{P_-} = 2^{a_1} R(-1)$ is divisible by 16. In the first case,
we get an isotropic element in $A_{P_+}[2]$ and then use Lemma~\ref{lem:isom quadratic spaces} to transfer it to $A_{P_-}[2]$; in the second case we use \cite[Proposition~1.8.2]{Nik}
to see directly that $A_{P_-}[2]$, being of order at least 16, must itself contain an isotropic vector.
\end{proof}


\begin{remark}
There are a few additional cases listed in Table~\ref{table:L-polys3} in which $P_-$ has rank 16 and admits an \emph{odd} overlattice of index $2$. In these cases we see that the overlattice represents neither 1 (the orthogonal complement of a norm-1 vector would have an even sublattice of discriminant large enough to violate Proposition~\ref{prop:rootless lattice discriminant} with $k=15$) nor 2 (because any root belongs to the even sublattice $P_-$). However, here again we are thwarted by the Barnes--Wall lattice, which also admits an odd overlattice of index 2 representing neither 1 nor 2 (the \emph{odd Barnes--Wall lattice}).

To outflank this adversary, we must adopt an approach not based solely on the size of the discriminant. One option would be the approach of \cite{Elkies} to root-counting using modular forms in place of linear programming, as suggested in Remark~\ref{remark:elkies bound}. We will instead consider embeddings into Niemeier lattices; see Remark~\ref{remark:niemeier embedding}.
\end{remark}

\section{Niemeier embeddings}

At this point we are left with a handful of cases which require further attention,
although the size of their discriminants suggest that it should still be impossible to find an even rootless lattice $P_-$ with the correct rank and discriminant.
A natural strategy for investigating this is to consider embeddings into Niemeier lattices, as in \cite{Kondo}.

\begin{lemma} \label{L:genus pminus}
    For each $L$-polynomial listed in Table~\ref{table:L-polys3}, the number of possible genera of the lattice $P_-$ is $2^{\Omega(R(-1))}$ where $\Omega(n)$ denotes the number of prime factors of the positive integer $n$ (counted with multiplicity).
    The corresponding discriminant forms are listed in Table~\ref{table:L-polys3a}.
\end{lemma}
\begin{proof}
Since $\rk(P_-) \geq 16$, the genus of $P_-$ is uniquely determined by the rank, signature, and discriminant form, and is composed of separate factors for each prime. 
When $R(-1)$ is exactly divisible by $p$ for some odd prime $p$,
there are two options for the genus at $p$, corresponding to whether or not the norm form on $A_K[p]$ is isotropic;
in the one case where $R(-1)$ is exactly divisible by $p^2$, we get 4 options.
In addition, the 2-part of the genus is constrained by 
Lemma~\ref{lem:isom quadratic spaces} and Lemma~\ref{lem;cyclic overlattice constraint}:
when $R(-1)$ is odd the genus is determined completely by $P_+$, and otherwise there are two possibilities. Finally there is a single parity condition constraining all of the local genera at once (the \emph{oddity condition} of \cite[Chapter 15, \S 7.7]{CS}). We deduce that there are $2^{\Omega(R(-1))}$ possible genera.

To check that the discriminant forms in Table~\ref{table:L-polys3a} are correct, we first observe directly that the $2$-adic components satisfy the constraints listed above.
We then verify that the lattice $K$ has matching discriminant form (see the proof of Corollary~\ref{cor:Niemeier embedding constraint}); this ensures that the oddity condition is satisfied.
\end{proof}

\begin{lemma} \label{lem:niemeier embedding}
Suppose that the $L$-polynomial of $X$ appears in Table~\ref{table:L-polys3}.
Then there exists an even positive definite lattice $K$ of rank $24 - a_2$ containing a root such that
any lattice in the genus of $K$ occurs as the orthogonal complement of $P_-(-1)$ for some primitive embedding into a Niemeier lattice $N$ (which cannot be the Leech lattice). 
\end{lemma}
\begin{proof}
In each case, the odd part of $R(-1)$ is squarefree, so $\ell(A_{P_-}) \leq 2 + v_2(R(-1))$; in particular, $\rk(P_-) + \ell(A_{P_-}) \leq 19$.
We may apply Nikulin's embedding theorem \cite[Corollary~1.12.3]{Nik} to conclude that $P_-(-1) \oplus A_1$ embeds primitively into some even unimodular lattice of rank 24, which must have a root because $A_1$ does. 
Let $K$ be the orthogonal complement of $P_-(-1)$ via this embedding, which by design contains $A_1$.
By Nikulin's theorem again,
any lattice in the same genus as $K$ is also the orthogonal complement of $P_-$ via a primitive embedding into some even unimodular lattice of rank $24$. This completes the proof.
\end{proof}

\begin{remark} \label{remark:niemeier embedding}
The trick of forcing a root into $N$ to exclude the Leech lattice occurs commonly in the study of automorphisms of K3 surfaces (e.g., \cite{Kondo}). In our case it has the important effect of ridding us of the ghost of the Barnes--Wall lattice $BW$: since $A_{BW} \cong (\ZZ/2\ZZ)^8$, we cannot apply Nikulin's theorem to $BW \oplus A_1$ to force an embedding of $BW$ into a Niemeier lattice with roots.

In the cases we are considering, we can force even more roots into $K$: 
the lattice $P_-(-1) \oplus D_4$ admits an even index-$2$ overlattice with $A_{P_-(-1) \oplus D_4} \cong A_{P_-(-1)}$, so it can also be embedded into a Niemeier lattice using Nikulin's theorem. However, doing this would be counterproductive for us, as we will need to optimize the choice of $K$ differently in order to apply Lemma~\ref{lem:niemeier LP setup}.
\end{remark}

With notation as in Lemma~\ref{lem:niemeier embedding},
at this point one possible approach would be to compute primitive embeddings of $K$ into all Niemeier lattices (e.g. using the functionality implemented in the \texttt{Oscar} computer algebra system), then directly observe that the orthogonal complement fails to be rootless in all cases.
We take a different approach based on linear programming duality, which requires remarkably little information about the Niemeier lattices themselves.

\begin{lemma} \label{lem:niemeier LP setup}
let $K$ be an even positive defiinte lattice with a root, and let $M_K$ be the Gram matrix of $K$ with respect to some basis.
Suppose that there exists a symmetric positive definite $k \times k$ integer matrix $T$ with even diagonal entries such that:
\begin{enumerate}
    \item[(a)] $\Trace(T M_K) \leq 24$; and
    \item[(b)] if $\Trace(T M_K) = 24$, then there exists $r \in R_K$ such that $\mathbf{v}_r^T T \mathbf{v}_r > 2$.
\end{enumerate}
Then there exists no primitive embedding of $K$ into a Niemeier lattice with rootless orthogonal complement.
\end{lemma}
\begin{proof}
Suppose to the contrary that $K \hookrightarrow N$ is a primitive embedding into a Niemeier lattice $N$ with rootless orthogonal complement. Since $K$ has a root, $N$ canot be the Leech lattice.
Let $h$ be the common Coxeter number of all indecomposable components of the root lattice $R(N)$.
Let $R_K, R_N$ be the sets of roots of $K,N$, respectively;
then 
\begin{equation} \label{eq:count roots}
\#R_N = h \rk N = 24 h.
\end{equation}
Let $e_1,\dots,e_k$ be the chosen basis of $K$.
For each $r \in R_N$, 
define the column vector 
\[
\mathbf{v}_r := (\langle e_1, r \rangle_K, \dots, \langle e_k, r\rangle_K) \in \ZZ^k.
\]
The Niemeier--Venkov identity\footnote{This can be interpreted as the statement that every Niemeier lattice is \emph{eutactic}.} \cite[Proposition 18.1]{CS} expresses the pairing on $N$ in terms of $R_N$; restricting it to $K$ yields
\begin{equation} \label{eq:venkov}
\sum_{r \in R_N} \mathbf{v}_r \mathbf{v}_r^T = 2hM_K.
\end{equation}
The condition that $K$ has rootless orthogonal complement
implies that $\mathbf{v}_r \neq 0$ for all $r \in R_N$.
Since $T$ is integral, symmetric, positive definite, and has even diagonal entries, it follows that $\mathbf{v}_r^T T \mathbf{v}_r$ is a positive even integer,
so $\mathbf{v}_r^T T \mathbf{v}_r - 2$ is a nonnegative integer.

On the other hand, by \eqref{eq:count roots}, \eqref{eq:venkov}, and condition (a), we have
\begin{align*}
\sum_{r \in R_N} (\mathbf{v}_r^T T \mathbf{v}_r - 2) &= 
\sum_{r \in R_N} (\Trace(\mathbf{v}_r^T T \mathbf{v}_r)-2) \\
&= \sum_{r \in R_N} (\Trace(T \mathbf{v}_r \mathbf{v}_r^T) - 2)\\
&= \Trace\left( T \sum_{r \in R_N} \mathbf{v}_r \mathbf{v}_r^T\right) - \sum_{r \in R_N} 2 \\
&= 2h \Trace(T M_K) - 2 \#R_N = 2h (\Trace(T M_K) - 24) \leq 0.
\end{align*}
This forces both $\Trace(T M_K) = 24$ and $\mathbf{v}_r^T T \mathbf{v}_r = 2$ for all $r \in R_N$ (and in particular for all $r \in R_K$), violating (b).
This yields the desired contradiction.
\end{proof}

\begin{remark}
When applying Lemma~\ref{lem:niemeier LP setup}, we have wide latitude in choosing the matrix $T$.
However, since $TM_K$ is again positive definite (e.g., because $T$ has a real square root), we can apply the arithmetic-geometric mean inequality to its eigenvalues; condition (a) therefore implies 
\[
\det(T) \det(M_K) = \det(TM_K) \leq \left( \frac{\Trace(TM_K)}{k} \right)^{k} \leq \left( \frac{24}{k} \right)^k.
\]
This suggests choosing $T$ to have as small determinant as possible; for $k = 1, \dots, 8$ this suggests taking $T$ to represent the root lattice $A_1, A_2, A_3, D_4, D_5, E_6, E_7, E_8$, respectively.
\end{remark}

\begin{lemma} \label{cor:Niemeier embedding constraint}
There is no K3 surface over a finite field with $L$-polynomial appearing in Table~\ref{table:L-polys3}.
(This concludes the proof of Theorem~\ref{thm:main}.)
\end{lemma}
\begin{proof}
We apply Lemma~\ref{lem:niemeier embedding} with $K$ as indicated in Table~\ref{table:L-polys3a}; $r = e_k$ equal to the last basis vector;
and $T$ the Gram matrix of the root lattice $E_7$ or $E_8$ for the basis corresponding to this labeling of the Dynkin diagrams:
\[
    \begin{tikzpicture}
        
        \node[circle,draw,label={[yshift=-1cm]1}] at (0,0) (t1) {};    
        \node[circle,draw,label={[yshift=-1cm]2}] at (2,0) (t2) {};    
        \node[circle,draw,label={[yshift=-1cm]3}] at (4,0) (t3) {};
        \node[circle,draw,label={[xshift=0.5cm,yshift=-0.5cm]4}] at (4,2) (t4) {};
        \node[circle,draw,label={[yshift=-1.1cm]5}] at (6,0) (t5) {};
        \node[circle,draw,label={[yshift=-1.1cm]6}] at (8,0) (t6) {};
        \node[circle,draw,label={[yshift=-1.1cm]7}] at (10,0) (t7) {};
        \node[circle,draw,label={[yshift=-1.1cm]8}] at (12,0) (t8) {};        
        
        \draw [-] (t1) edge (t2) (t2) edge (t3) (t3) edge (t4) (t3) edge (t5) (t5) edge (t6) (t6) edge (t7) (t7) edge [dashed] (t8); 
        
    \end{tikzpicture}
\]
A verification script is provided in the GitHub repository to confirm that each lattice $K$ satisfies $A_K \cong A_{P_-}$, $\Trace(TM_K) = 24$, and $\mathbf{v}_r^T T \mathbf{v}_T > 2$.
\end{proof}

\begin{remark} \label{rem:bounds from Niemeier embeddings}
One could also use Lemma~\ref{lem:niemeier LP setup} to directly improve the bounds in Proposition~\ref{prop:rootless lattice discriminant}. While this would be beneficial for further study of K3 surfaces, we will not take up this task here.
\end{remark}

\section{Construction of pointless K3 surfaces}

We conclude by verifying the positive assertion of Theorem~\ref{thm:main} for $q=2,3,4,5,7,8,9,11$. For this one can give a uniform construction whose shape is dictated by Proposition~\ref{prop:square 0 or 2}.

\begin{proposition} \label{prop:pointless k3 over F11}
    For $q \leq 11$, there exists a K3 surface $X$ over $\FF_q$ admitting an elliptic fibration satisfying the following conditions.
    \begin{itemize}
        \item[(a)] The fiber over each point of $\PP^1(\FF_q)$ is singular of type $I_2$.
        (For $q = 11$ this accounts for all singular fibers.)
        \item[(b)] No such fiber has a point over $\FF_q$.
    \end{itemize}
    Consequently, $X$ is pointless.
\end{proposition}
\begin{proof}
    We construct $X$ as a smooth double cover of $Y := \PP^1_{\FF_q} \times_{\FF_q} \PP^1_{\FF_q} = \Proj \FF_q[X,Y] \times_{\FF_q} \Proj {\FF_q}[Z,W]$, using the fibration over the $Z,W$-line. This double cover will be the zero locus of 
    $T^2 + T H = F$
    for some bihomogeneous forms $H$ of bidegree $(2,2)$, taken to be 0 if $q$ is odd, and $F$ of bidegree $(4,4)$. To achieve conditions (a) and (b),
    we will arrange that for each $t \in \mathbb{P}^1(\FF_q)$:
    \begin{itemize}
        \item if $q$ is odd, then there exists a separable irreducible binary quadratic form $Q_t \in \FF_q[X, Y]$ such that $F_t = c_t Q_t^2$ for some nonsquare $c_t \in \FF_q^\times$;
        \item if $q$ is even, then $H_t$ is itself an irreducible binary quadratic form and $F_t = R_t^2 + R_t H_t + c_t H_t^2$ for some quadratic form $R_t$ and some $c_t \in \FF_q$ of absolute trace 1.
    \end{itemize}
    The fiber at $t$ will then be geometrically of type $I_2$, consisting of two copies of $\PP^1$ meeting transversally at two points (sometimes called a ``banana''), where neither the curves nor the crossing points are fixed by Frobenius (i.e., a ``banana nonsplit'').
    We provide equations in Table~\ref{table:pointless eqs1} (for $q$ odd) and Table~\ref{table:pointless eqs2} (for $q$ even) meeting these conditions;
    a verification script using Magma can be found in our GitHub repository.
\end{proof}

\begin{remark} \label{rem:pointless K3 over F11}
    We elaborate on Proposition~\ref{prop:pointless k3 over F11} in the cases $q = 3,5,7,11$
    (the case $q=9$ requires a somewhat different setup, as do the characteristic 2 cases).
    We take
    \[
      F = Q_2^2 R_2 - Q_1^2 R_1
    \]
    where $Q_1, Q_2 \in \FF_q[X,Y]$ are irreducible binary quadratic forms and
    $R_1, R_2 \in \FF_q[Z,W]$ are nonsingular quartic forms. 
    This setup let us take care of the fiber condition at many points at once:
    \begin{itemize}
        \item if $R_1$ vanishes at $t$, then $R_2$ must specialize to a nonsquare at $t$;
        \item if $R_2$ vanishes at $t$, then $-R_1$ must specialize to a nonsquare at $t$.
    \end{itemize}
    This suffices to handle $q = 3,5,7$ as we can ensure that $R_1R_2$ vanishes at all points of $\PP^1(\FF_q)$. 
    
    For $q = 11$, we can ensure that $R_1 R_2 (R_1+R_2)$ vanishes at all points of $\PP^1(\FF_q)$.
    To cover the points where $R_1+R_2$ vanishes, we arrange that $Q_1^2 + Q_2^2 = 2Q_3^2$ for some third irreducible binary quadratic form $Q_3$.
    Then our defining equation can be rewritten as
    \[
    T^2 = Q_2^2(R_1+R_2) - 2Q_3^2 R_1;
    \]
    if $R_1+R_2$ vanishes at $t$, then $-2R_1$ must specialize to a nonsquare at $t$. 
\end{remark}

\begin{table}
    \label{table:pointless eqs1}
    \caption{Equations $T^2 = F$ for pointless K3 surfaces over $\FF_q$ for odd $q$. For $q = 9$, $i \in \FF_9$ satisfies $i^2 = -1$ and $d:=i-1$.}
    {\scriptsize
    \begin{tabular}{c|cccc}
        $q$ & $Q_1$ & $Q_2$ & $R_1$ & $F - R_1 Q_1^2$ \\
        \hline
        $3$ & $X^2+Y^2$ & $X^2+XY-Y^2$  & $2(Z+W)(Z^3-Z^2W+W^3)$ & $2ZW(Z^2+W^2)Q_2^2$\\
        $5$ & $X^2+2Y^2$ & $X^2+XY+Y^2$ & $2(Z-3W)(Z-4W)(Z^2+3W^2)$ & $3ZW(Z-W)(Z-2W)Q_2^2$ \\
        $7$ & $X^2+Y^2$ & $2X^2-XY-2Y^2$ & $3W(Z-3W)(Z-5W)(Z-6W)$ & $6Z(Z-W)(Z-2W)(Z-4W)Q_2^2$ \\
        $11$ & $X^2 -2XY-Y^2$ & $X^2+2XY-Y^2$  & $6(Z-3W)(Z-4W)(Z-6W)(Z-8W)$ & $5WZ(Z-W)(Z-2W)Q_2^2$ \\
        \hline
        $9$ & $X^2 + dY^2$ & $X^2 - dY^2$ & $(1+i)(W^4+i Z^4)$ & 
        $ZWQ_2(ZWQ_2+ i (W^2-dZ^2)XY)$ \\
    \end{tabular}}
    
\end{table}

\begin{table}
    \label{table:pointless eqs2}
    \caption{Equations $T^2+HT = F$ for pointless K3 surfaces over $\FF_q$ for even $q$. We specify the bihomogeneous forms $H$ and $F$ as matrices.  For $q = 4,8$, $\alpha \in \FF_4$ satisfies $\alpha^2 = \alpha+1$ and $\beta \in \FF_8$ satisfies $\beta^3 = \beta+1$.}

    \begin{tabular}{c|cc}
     $q$ & $H$ & $F$ \\
     \hline
     $2$ & $\begin{pmatrix}
         1 & 1 & 1 \\ 1 & 1 & 1 \\ 1 & 1 & 1
     \end{pmatrix}$  & \scriptsize $\begin{pmatrix}
         1 & 0 & 1 & 0 & 1 \\
         0 & 0 & 1 & 1 & 0 \\
         1 & 1 & 1 & 1 & 1\\
         0 & 1 & 1 & 0 & 0\\
         1 & 0 & 1 & 0 & 1
     \end{pmatrix}$ \\
     $4$ & $\begin{pmatrix} 1+\alpha & 1 & 1 \\ 1 & 1 & 1+\alpha \\ 1 & \alpha & 1 \end{pmatrix}$ & \scriptsize
     $\begin{pmatrix} 1 & 1+\alpha & 0 & 0 & 1+\alpha \\ \alpha & 1+\alpha & 0 & 0 & 1+\alpha \\ 0 & 0 & 0 & 0 & 0 \\ 0 & 0 & 0 & 0 & 0 \\
     \alpha & 1+\alpha & 0 & 0 & 1+\alpha \end{pmatrix}$\\
     $8$ & 
     $\begin{pmatrix}
         \beta^2 & \beta^2 & 1+\beta+\beta^2 \\
         \beta^2& \beta^2 & 1 \\
         \beta & \beta^2 & \beta
     \end{pmatrix}$ & \scriptsize $\begin{pmatrix} 1 & 1 & 0 & 0 & 1+\beta^2 \\ 1 & 1 & 0 & 0 & 1+\beta+\beta^2 \\ 0 & 0 & 0 & 0 & 0 \\
     0 & 0 & 0 & 0 & 0 \\ 1 & 1 & 0 & 0 & 1 \end{pmatrix}$
    \end{tabular}
\end{table}

\bibliographystyle{amsplain}
\bibliography{pointless_bibliography}

\begin{table}[ht]
\caption{Certificates for Lemma~\ref{lem:niemeier LP setup} for genera indicated in Table~\ref{table:L-polys3}. Each Gram matrix is labeled above by the Gram matrix of $P_+$ and the discriminant form $(A_{P_-}, q_{P_-})$ expressed as a direct sum: $\tfrac{a}{n}$ represents the form $\ZZ/n\ZZ \to \QQ/2\ZZ$ taking $1$ to $\tfrac{a}{n}$, while $\mathcal{U}$ represents the form $(\ZZ/2\ZZ)^2 \to \QQ/2\ZZ$
with upper triangular Gram matrix $(1,1/2,1)$. We write $-$ in place of $-1$.}
\label{table:L-polys3a}

\renewcommand*{\arraystretch}{2}
\begin{tabular}{cccc}
\setlength\arraycolsep{2pt}
\renewcommand*{\arraystretch}{0.75}
$\scriptsize\begin{pmatrix} 6 \end{pmatrix}$ $(\tfrac{1}{2}, \tfrac{8}{5}, \tfrac{60}{31})$
& 
\setlength\arraycolsep{2pt}
\renewcommand*{\arraystretch}{0.75}
$\scriptsize\begin{pmatrix} 6 \end{pmatrix}$ $(\tfrac{1}{2}, \tfrac{6}{5}, \tfrac{56}{31})$ &
\setlength\arraycolsep{2pt}
\renewcommand*{\arraystretch}{0.75}
$\scriptsize\begin{pmatrix} 4 \end{pmatrix}$ $(\tfrac{3}{4}, \tfrac{12}{7}, \tfrac{20}{11})$
&
\setlength\arraycolsep{2pt}
\renewcommand*{\arraystretch}{0.75}
$\scriptsize\begin{pmatrix} 4 \end{pmatrix}$ $(\tfrac{3}{4}, \tfrac{8}{7}, \tfrac{18}{11})$
\\
\setlength\arraycolsep{2pt}
\renewcommand*{\arraystretch}{0.75}
$\scriptsize \begin{pmatrix}
2 & 2 & 2 & 1 & 1 & 1 & 1 \\
2 & 4 & 4 & 2 & 3 & 2 & 1 \\
2 & 4 & 10 & 4 & 6 & 3 & 1 \\
1 & 2 & 4 & 4 & 3 & 2 & 1 \\
1 & 3 & 6 & 3 & 6 & 3 & 1 \\
1 & 2 & 3 & 2 & 3 & 4 & 1 \\
1 & 1 & 1 & 1 & 1 & 1 & 2
\end{pmatrix}$ &
\setlength\arraycolsep{2pt}
\renewcommand*{\arraystretch}{0.75}
$\scriptsize \begin{pmatrix}
2 & 2 & 2 & 1 & 1 & 1 & 1 \\
2 & 6 & 4 & 2 & 2 & 1 & 1 \\
2 & 4 & 8 & 4 & 5 & 2 & 1 \\
1 & 2 & 4 & 4 & 2 & 0 & 0 \\
1 & 2 & 5 & 2 & 6 & 3 & 1 \\
1 & 1 & 2 & 0 & 3 & 4 & 2 \\
1 & 1 & 1 & 0 & 1 & 2 & 2
\end{pmatrix}$ &
\setlength\arraycolsep{2pt}
\renewcommand*{\arraystretch}{0.75}
$\scriptsize \begin{pmatrix}
2 & 2 & 2 & 1 & 1 & 1 & 0 \\
2 & 6 & 6 & 3 & 3 & 2 & 1 \\
2 & 6 & 10 & 5 & 5 & 4 & 2 \\
1 & 3 & 5 & 4 & 2 & 2 & 1 \\
1 & 3 & 5 & 2 & 8 & 4 & 2 \\
1 & 2 & 4 & 2 & 4 & 4 & 2 \\
0 & 1 & 2 & 1 & 2 & 2 & 2
\end{pmatrix}$ &
\setlength\arraycolsep{2pt}
\renewcommand*{\arraystretch}{0.75}
$\scriptsize \begin{pmatrix}
4 & 4 & 3 & 2 & 1 & 1 & - \\
4 & 8 & 8 & 4 & 4 & 2 & - \\
3 & 8 & 12 & 6 & 7 & 4 & 0 \\
2 & 4 & 6 & 4 & 4 & 2 & 0 \\
1 & 4 & 7 & 4 & 8 & 4 & 1 \\
1 & 2 & 4 & 2 & 4 & 4 & 1 \\
- & - & 0 & 0 & 1 & 1 & 2
\end{pmatrix}$ 
\\
\setlength\arraycolsep{2pt}
\renewcommand*{\arraystretch}{0.75}
$\scriptsize\begin{pmatrix} 4 \end{pmatrix}$ $(\tfrac{7}{4}, \tfrac{12}{7}, \tfrac{18}{11})$
&
\setlength\arraycolsep{2pt}
\renewcommand*{\arraystretch}{0.75}
$\scriptsize\begin{pmatrix} 4 \end{pmatrix}$ $(\tfrac{7}{4}, \tfrac{8}{7}, \tfrac{20}{11})$ &
\setlength\arraycolsep{2pt}
\renewcommand*{\arraystretch}{0.75}
$\scriptsize\begin{pmatrix} 4 & 0 \\ 0 & -30 \end{pmatrix}$ $(\tfrac{5}{4}, \tfrac{1}{2}, \tfrac{4}{3}, \tfrac{22}{13})$
& 
\setlength\arraycolsep{2pt}
\renewcommand*{\arraystretch}{0.75}
$\scriptsize\begin{pmatrix} 4 & 0 \\ 0 & -30 \end{pmatrix}$ $(\tfrac{5}{4}, \tfrac{1}{2}, \tfrac{2}{3}, \tfrac{24}{13})$ 
\\
\setlength\arraycolsep{2pt}
\renewcommand*{\arraystretch}{0.75}
$\scriptsize \begin{pmatrix}
4 & 3 & 2 & 0 & 1 & 1 & 0 \\
3 & 6 & 5 & 2 & 2 & 1 & 0 \\
2 & 5 & 8 & 4 & 4 & 1 & 0 \\
0 & 2 & 4 & 4 & 2 & 0 & 0 \\
1 & 2 & 4 & 2 & 4 & 1 & 1 \\
1 & 1 & 1 & 0 & 1 & 2 & 1 \\
0 & 0 & 0 & 0 & 1 & 1 & 2
\end{pmatrix}$ &
\setlength\arraycolsep{2pt}
\renewcommand*{\arraystretch}{0.75}
$\scriptsize \begin{pmatrix}
2 & 1 & 2 & 1 & 1 & 0 & 0 \\
1 & 6 & 5 & 2 & 3 & 2 & 1 \\
2 & 5 & 10 & 4 & 6 & 4 & 2 \\
1 & 2 & 4 & 4 & 3 & 2 & 1 \\
1 & 3 & 6 & 3 & 6 & 4 & 2 \\
0 & 2 & 4 & 2 & 4 & 4 & 2 \\
0 & 1 & 2 & 1 & 2 & 2 & 2
\end{pmatrix}$ &
\setlength\arraycolsep{2pt}
\renewcommand*{\arraystretch}{0.75}
$\scriptsize \begin{pmatrix}
6 & 4 & 5 & 2 & 3 & 2 & 1 & 0 \\
4 & 8 & 10 & 5 & 7 & 5 & 3 & 1 \\
5 & 10 & 16 & 8 & 11 & 7 & 4 & 2 \\
2 & 5 & 8 & 6 & 5 & 3 & 2 & 1 \\
3 & 7 & 11 & 5 & 10 & 7 & 4 & 2 \\
2 & 5 & 7 & 3 & 7 & 6 & 4 & 2 \\
1 & 3 & 4 & 2 & 4 & 4 & 4 & 2 \\
0 & 1 & 2 & 1 & 2 & 2 & 2 & 2
\end{pmatrix}$ &
\setlength\arraycolsep{2pt}
\renewcommand*{\arraystretch}{0.75}
$\scriptsize \begin{pmatrix}
4 & 4 & 4 & 1 & 2 & 1 & 1 & 0 \\
4 & 8 & 8 & 4 & 5 & 2 & 1 & 0 \\
4 & 8 & 12 & 5 & 8 & 4 & 2 & 1 \\
1 & 4 & 5 & 4 & 4 & 2 & 1 & 0 \\
2 & 5 & 8 & 4 & 8 & 4 & 2 & 1 \\
1 & 2 & 4 & 2 & 4 & 4 & 2 & 1 \\
1 & 1 & 2 & 1 & 2 & 2 & 2 & 1 \\
0 & 0 & 1 & 0 & 1 & 1 & 1 & 2
\end{pmatrix}$ \\

\setlength\arraycolsep{2pt}
\renewcommand*{\arraystretch}{0.75}
$\scriptsize\begin{pmatrix} 4 & 0 \\ 0 & -30 \end{pmatrix}$ $(\tfrac{1}{4}, \tfrac{1}{2}, \tfrac{4}{3}, \tfrac{24}{13})$ 
&
\setlength\arraycolsep{2pt}
\renewcommand*{\arraystretch}{0.75}
$\scriptsize\begin{pmatrix} 4 & 0 \\ 0 & -30 \end{pmatrix}$ $(\tfrac{1}{4}, \tfrac{1}{2}, \tfrac{2}{3}, \tfrac{22}{13})$ &
\setlength\arraycolsep{2pt}
\renewcommand*{\arraystretch}{0.75}
$\scriptsize\begin{pmatrix} 4 & 2 \\ 2 & -28 \end{pmatrix}$  $(\mathcal{U}, \tfrac{12}{7}, \tfrac{20}{11})$
& 
\setlength\arraycolsep{2pt}
\renewcommand*{\arraystretch}{0.75}
$\scriptsize\begin{pmatrix} 4 & 2 \\ 2 & -28 \end{pmatrix}$  $(\mathcal{U}, \tfrac{8}{7}, \tfrac{18}{11})$ \\
\setlength\arraycolsep{2pt}
\renewcommand*{\arraystretch}{0.75}
$\scriptsize \begin{pmatrix}
4 & 4 & 4 & 2 & 3 & 2 & 2 & 1 \\
4 & 8 & 9 & 4 & 7 & 5 & 3 & 2 \\
4 & 9 & 14 & 6 & 11 & 8 & 5 & 3 \\
2 & 4 & 6 & 4 & 5 & 4 & 2 & 1 \\
3 & 7 & 11 & 5 & 10 & 7 & 4 & 3 \\
2 & 5 & 8 & 4 & 7 & 8 & 5 & 3 \\
2 & 3 & 5 & 2 & 4 & 5 & 6 & 2 \\
1 & 2 & 3 & 1 & 3 & 3 & 2 & 2
\end{pmatrix}$ &
\setlength\arraycolsep{2pt}
\renewcommand*{\arraystretch}{0.75}
$\scriptsize \begin{pmatrix}
2 & 3 & 3 & 1 & 3 & 2 & 1 & 0 \\
3 & 8 & 8 & 4 & 6 & 4 & 2 & 0 \\
3 & 8 & 14 & 6 & 10 & 6 & 3 & 1 \\
1 & 4 & 6 & 4 & 4 & 3 & 2 & 1 \\
3 & 6 & 10 & 4 & 10 & 6 & 3 & 1 \\
2 & 4 & 6 & 3 & 6 & 6 & 3 & 1 \\
1 & 2 & 3 & 2 & 3 & 3 & 4 & 2 \\
0 & 0 & 1 & 1 & 1 & 1 & 2 & 2
\end{pmatrix}$ &
\setlength\arraycolsep{2pt}
\renewcommand*{\arraystretch}{0.75}
$\scriptsize \begin{pmatrix}
4 & 5 & 6 & 2 & 4 & 2 & 1 & 0 \\
5 & 10 & 12 & 5 & 9 & 5 & 3 & 1 \\
6 & 12 & 18 & 9 & 13 & 8 & 5 & 2 \\
2 & 5 & 9 & 8 & 6 & 4 & 3 & 1 \\
4 & 9 & 13 & 6 & 12 & 8 & 5 & 2 \\
2 & 5 & 8 & 4 & 8 & 8 & 5 & 3 \\
1 & 3 & 5 & 3 & 5 & 5 & 4 & 2 \\
0 & 1 & 2 & 1 & 2 & 3 & 2 & 2
\end{pmatrix}$ &
\setlength\arraycolsep{2pt}
\renewcommand*{\arraystretch}{0.75}
$\scriptsize \begin{pmatrix}
4 & 5 & 6 & 3 & 4 & 2 & 1 & 0 \\
5 & 10 & 13 & 6 & 8 & 5 & 3 & 1 \\
6 & 13 & 20 & 9 & 12 & 8 & 5 & 2 \\
3 & 6 & 9 & 6 & 6 & 5 & 3 & 1 \\
4 & 8 & 12 & 6 & 12 & 8 & 5 & 3 \\
2 & 5 & 8 & 5 & 8 & 8 & 5 & 3 \\
1 & 3 & 5 & 3 & 5 & 5 & 4 & 2 \\
0 & 1 & 2 & 1 & 3 & 3 & 2 & 2
\end{pmatrix}$ \\
\setlength\arraycolsep{2pt}
\renewcommand*{\arraystretch}{0.75}
$\scriptsize\begin{pmatrix} 6 & 6 \\ 6 & -12 \end{pmatrix}$ $(\tfrac{1}{2}, \tfrac{1}{2}, \tfrac{4}{3}, \tfrac{38}{25})$
& 
\setlength\arraycolsep{2pt}
\renewcommand*{\arraystretch}{0.75}
$\scriptsize\begin{pmatrix} 6 & 6 \\ 6 & -12 \end{pmatrix}$ $(\tfrac{1}{2}, \tfrac{1}{2}, \tfrac{4}{3}, \tfrac{46}{25})$
&
\setlength\arraycolsep{2pt}
\renewcommand*{\arraystretch}{0.75}
$\scriptsize\begin{pmatrix} 6 & 6 \\ 6 & -12 \end{pmatrix}$ $(\tfrac{1}{2}, \tfrac{1}{2}, \tfrac{4}{3}, \tfrac{2}{5}, \tfrac{2}{5})$
& 
\setlength\arraycolsep{2pt}
\renewcommand*{\arraystretch}{0.75}
$\scriptsize\begin{pmatrix} 6 & 6 \\ 6 & -12 \end{pmatrix}$ $(\tfrac{1}{2}, \tfrac{1}{2}, \tfrac{4}{3}, \tfrac{2}{5}, \tfrac{2}{5})$ \\
\setlength\arraycolsep{2pt}
\renewcommand*{\arraystretch}{0.75}
$\scriptsize \begin{pmatrix}
4 & 4 & 4 & 2 & 3 & 2 & 1 & 0 \\
4 & 8 & 8 & 3 & 5 & 3 & 2 & 0 \\
4 & 8 & 10 & 4 & 6 & 4 & 2 & 0 \\
2 & 3 & 4 & 4 & 3 & 2 & 1 & 0 \\
3 & 5 & 6 & 3 & 6 & 4 & 1 & 0 \\
2 & 3 & 4 & 2 & 4 & 4 & 1 & 0 \\
1 & 2 & 2 & 1 & 1 & 1 & 2 & 1 \\
0 & 0 & 0 & 0 & 0 & 0 & 1 & 2
\end{pmatrix}$ &
\setlength\arraycolsep{2pt}
\renewcommand*{\arraystretch}{0.75}
$\scriptsize \begin{pmatrix}
4 & 4 & 4 & 2 & 3 & 2 & 2 & 0 \\
4 & 8 & 8 & 4 & 5 & 3 & 2 & 1 \\
4 & 8 & 10 & 5 & 6 & 4 & 2 & 1 \\
2 & 4 & 5 & 4 & 3 & 2 & 1 & 0 \\
3 & 5 & 6 & 3 & 6 & 4 & 3 & 1 \\
2 & 3 & 4 & 2 & 4 & 4 & 2 & 1 \\
2 & 2 & 2 & 1 & 3 & 2 & 4 & 1 \\
0 & 1 & 1 & 0 & 1 & 1 & 1 & 2
\end{pmatrix}$ &
\setlength\arraycolsep{2pt}
\renewcommand*{\arraystretch}{0.75}
$\scriptsize \begin{pmatrix}
2 & 2 & 2 & 1 & 2 & 1 & 0 & 0 \\
2 & 6 & 7 & 3 & 5 & 3 & 2 & 1 \\
2 & 7 & 12 & 5 & 8 & 6 & 4 & 2 \\
1 & 3 & 5 & 6 & 4 & 3 & 2 & 1 \\
2 & 5 & 8 & 4 & 8 & 6 & 4 & 2 \\
1 & 3 & 6 & 3 & 6 & 6 & 4 & 2 \\
0 & 2 & 4 & 2 & 4 & 4 & 4 & 2 \\
0 & 1 & 2 & 1 & 2 & 2 & 2 & 2
\end{pmatrix}$ &
\setlength\arraycolsep{2pt}
\renewcommand*{\arraystretch}{0.75}
$\scriptsize \begin{pmatrix}
2 & 3 & 4 & 2 & 3 & 3 & 2 & 1 \\
3 & 8 & 10 & 5 & 7 & 6 & 3 & 1 \\
4 & 10 & 16 & 8 & 12 & 9 & 5 & 2 \\
2 & 5 & 8 & 6 & 5 & 4 & 2 & 1 \\
3 & 7 & 12 & 5 & 12 & 8 & 5 & 2 \\
3 & 6 & 9 & 4 & 8 & 8 & 4 & 2 \\
2 & 3 & 5 & 2 & 5 & 4 & 4 & 1 \\
1 & 1 & 2 & 1 & 2 & 2 & 1 & 2
\end{pmatrix}$ \\
\setlength\arraycolsep{2pt}
\renewcommand*{\arraystretch}{0.75}
$\scriptsize \begin{pmatrix} 4 & 0 \\ 0 & -26 \end{pmatrix}$ $(\tfrac{1}{4}, \tfrac{3}{2}, \tfrac{72}{37})$ &
\setlength\arraycolsep{2pt}
\renewcommand*{\arraystretch}{0.75}
$\scriptsize \begin{pmatrix} 4 & 0 \\ 0 & -26 \end{pmatrix}$ $(\tfrac{5}{4}, \tfrac{3}{2}, \tfrac{70}{37})$ & 
\setlength\arraycolsep{2pt}
\renewcommand*{\arraystretch}{0.75}
$\scriptsize \begin{pmatrix} 4 & 2 \\ 2 & -20 \end{pmatrix}$
$(\mathcal{U}, \tfrac{4}{3}, \tfrac{44}{23})$ &
\setlength\arraycolsep{2pt}
\renewcommand*{\arraystretch}{0.75}
$\scriptsize \begin{pmatrix} 4 & 2 \\ 2 & -20 \end{pmatrix}$ 
$(\mathcal{U}, \tfrac{2}{3}, \tfrac{36}{23})$ \\
\setlength\arraycolsep{2pt}
\renewcommand*{\arraystretch}{0.75}
$\scriptsize \begin{pmatrix}
2 & 3 & 4 & 2 & 3 & 3 & 2 & 1 \\
3 & 8 & 10 & 5 & 7 & 6 & 3 & 1 \\
4 & 10 & 16 & 8 & 12 & 9 & 5 & 2 \\
2 & 5 & 8 & 6 & 5 & 4 & 2 & 1 \\
3 & 7 & 12 & 5 & 12 & 8 & 5 & 2 \\
3 & 6 & 9 & 4 & 8 & 8 & 4 & 2 \\
2 & 3 & 5 & 2 & 5 & 4 & 4 & 1 \\
1 & 1 & 2 & 1 & 2 & 2 & 1 & 2
\end{pmatrix}$ &
\setlength\arraycolsep{2pt}
\renewcommand*{\arraystretch}{0.75}
$\scriptsize \begin{pmatrix}
6 & 5 & 6 & 3 & 4 & 2 & 1 & 0 \\
5 & 10 & 13 & 6 & 10 & 7 & 4 & 1 \\
6 & 13 & 20 & 9 & 14 & 10 & 6 & 2 \\
3 & 6 & 9 & 6 & 7 & 5 & 3 & 1 \\
4 & 10 & 14 & 7 & 12 & 9 & 5 & 2 \\
2 & 7 & 10 & 5 & 9 & 8 & 4 & 2 \\
1 & 4 & 6 & 3 & 5 & 4 & 4 & 2 \\
0 & 1 & 2 & 1 & 2 & 2 & 2 & 2
\end{pmatrix}$ &
\setlength\arraycolsep{2pt}
\renewcommand*{\arraystretch}{0.75}
$\scriptsize \begin{pmatrix}
4 & 4 & 4 & 2 & 2 & 1 & - & - \\
4 & 8 & 8 & 4 & 5 & 3 & 0 & 0 \\
4 & 8 & 12 & 6 & 8 & 4 & 2 & 0 \\
2 & 4 & 6 & 4 & 4 & 2 & 1 & 0 \\
2 & 5 & 8 & 4 & 8 & 5 & 3 & 1 \\
1 & 3 & 4 & 2 & 5 & 6 & 3 & 2 \\
- & 0 & 2 & 1 & 3 & 3 & 4 & 2 \\
- & 0 & 0 & 0 & 1 & 2 & 2 & 2
\end{pmatrix}$ &
\setlength\arraycolsep{2pt}
\renewcommand*{\arraystretch}{0.75}
$\scriptsize \begin{pmatrix}
4 & 4 & 4 & 2 & 2 & 2 & 0 & 1 \\
4 & 8 & 8 & 4 & 5 & 4 & 2 & 1 \\
4 & 8 & 12 & 6 & 9 & 7 & 4 & 2 \\
2 & 4 & 6 & 4 & 4 & 3 & 2 & 1 \\
2 & 5 & 9 & 4 & 10 & 7 & 4 & 2 \\
2 & 4 & 7 & 3 & 7 & 6 & 3 & 2 \\
0 & 2 & 4 & 2 & 4 & 3 & 4 & 1 \\
1 & 1 & 2 & 1 & 2 & 2 & 1 & 2
\end{pmatrix}$ \\
\setlength\arraycolsep{2pt}
\renewcommand*{\arraystretch}{0.75}
$\scriptsize \begin{pmatrix} 6 & 0 \\ 0 & -12 \end{pmatrix}$ $(\tfrac{5}{4}, \tfrac{3}{2}, \tfrac{2}{3}, \tfrac{18}{11})$ &
\setlength\arraycolsep{2pt}
\renewcommand*{\arraystretch}{0.75}
$\scriptsize \begin{pmatrix} 6 & 0 \\ 0 & -12 \end{pmatrix}$ $(\tfrac{5}{4}, \tfrac{3}{2}, \tfrac{4}{3}, \tfrac{20}{11})$ &
\setlength\arraycolsep{2pt}
\renewcommand*{\arraystretch}{0.75}
$\scriptsize \begin{pmatrix} 6 & 0 \\ 0 & -12 \end{pmatrix}$ $(\tfrac{1}{4}, \tfrac{3}{2},\tfrac{2}{3}, \tfrac{20}{11})$ &
\setlength\arraycolsep{2pt}
\renewcommand*{\arraystretch}{0.75}
$\scriptsize \begin{pmatrix} 6 & 0 \\ 0 & -12 \end{pmatrix}$ $( \tfrac{1}{4}, \tfrac{3}{2}, \tfrac{4}{3},\tfrac{18}{11})$ \\
\setlength\arraycolsep{2pt}
\renewcommand*{\arraystretch}{0.75}
$\scriptsize \begin{pmatrix}
2 & 2 & 3 & 2 & 2 & 1 & 0 & 0 \\
2 & 6 & 5 & 3 & 3 & 2 & 1 & 0 \\
3 & 5 & 10 & 5 & 7 & 5 & 3 & 1 \\
2 & 3 & 5 & 4 & 3 & 2 & 1 & 0 \\
2 & 3 & 7 & 3 & 8 & 6 & 3 & 1 \\
1 & 2 & 5 & 2 & 6 & 6 & 3 & 1 \\
0 & 1 & 3 & 1 & 3 & 3 & 4 & 2 \\
0 & 0 & 1 & 0 & 1 & 1 & 2 & 2
\end{pmatrix}$ &
\setlength\arraycolsep{2pt}
\renewcommand*{\arraystretch}{0.75}
$\scriptsize \begin{pmatrix}
2 & 3 & 4 & 2 & 3 & 2 & 2 & 1 \\
3 & 8 & 10 & 4 & 7 & 5 & 4 & 2 \\
4 & 10 & 16 & 8 & 12 & 9 & 6 & 3 \\
2 & 4 & 8 & 6 & 6 & 5 & 3 & 2 \\
3 & 7 & 12 & 6 & 14 & 9 & 5 & 2 \\
2 & 5 & 9 & 5 & 9 & 8 & 4 & 2 \\
2 & 4 & 6 & 3 & 5 & 4 & 4 & 2 \\
1 & 2 & 3 & 2 & 2 & 2 & 2 & 2
\end{pmatrix}$ &
\setlength\arraycolsep{2pt}
\renewcommand*{\arraystretch}{0.75}
$\scriptsize \begin{pmatrix}
2 & 2 & 2 & 1 & 1 & 0 & 0 & 0 \\
2 & 6 & 5 & 2 & 3 & 0 & 0 & 0 \\
2 & 5 & 8 & 4 & 5 & 2 & 1 & 0 \\
1 & 2 & 4 & 4 & 2 & 1 & 0 & 0 \\
1 & 3 & 5 & 2 & 6 & 3 & 1 & 0 \\
0 & 0 & 2 & 1 & 3 & 4 & 2 & 1 \\
0 & 0 & 1 & 0 & 1 & 2 & 2 & 1 \\
0 & 0 & 0 & 0 & 0 & 1 & 1 & 2
\end{pmatrix}$ &
\setlength\arraycolsep{2pt}
\renewcommand*{\arraystretch}{0.75}
$\scriptsize \begin{pmatrix}
4 & 3 & 3 & 2 & 2 & 1 & 0 & 0 \\
3 & 8 & 9 & 5 & 7 & 5 & 3 & 2 \\
3 & 9 & 16 & 8 & 11 & 8 & 5 & 2 \\
2 & 5 & 8 & 6 & 5 & 4 & 2 & 1 \\
2 & 7 & 11 & 5 & 10 & 7 & 5 & 2 \\
1 & 5 & 8 & 4 & 7 & 6 & 4 & 2 \\
0 & 3 & 5 & 2 & 5 & 4 & 4 & 2 \\
0 & 2 & 2 & 1 & 2 & 2 & 2 & 2
\end{pmatrix}$ \\
\end{tabular}
\end{table}

\end{document}